\documentclass[12pt,a4paper]{article}
\usepackage{fullpage,amsmath,amssymb,diagram,amsthm,verbatim,alltt,amsfonts,array}
\usepackage[english]{babel}

 \DeclareMathOperator{\Pf}{pf}
 \DeclareMathOperator{\Co}{Coker}
\DeclareMathOperator{\coker}{coker}
 \DeclareMathOperator{\Rk}{rank}
 \DeclareMathOperator{\sgn}{sign}
 \DeclareMathOperator{\im}{Im}
 \DeclareMathOperator{\codim}{codim}
 \DeclareMathOperator{\depth}{depth}
 \DeclareMathOperator{\id}{Id}
 \newcommand{\To}{\longrightarrow}
 \newcommand{\ra}{\rightarrow}
 \newcommand{\km}{\mathcal{M}}
 
\renewcommand{\P}{\mathbb{P}}
\newcommand{\A}{\mathbb{A}}

\date{}
\begin{document}

\theoremstyle{plain}
 \newtheorem{thm}{Theorem}[section]
 \newtheorem{cor}[thm]{Corollary}
 \newtheorem{lem}[thm]{Lemma}
 \newtheorem{prop}[thm]{Proposition}
 \newtheorem{defn}[thm]{Definition}

\theoremstyle{remark}
 \newtheorem{rem}[thm]{Remark}
 \newtheorem{rems}[thm]{Remarks}

\numberwithin{equation}{section}

\title{\textbf{Rank two Cohen-Macaulay modules over singularities of type
    $\mathbf{x_1^3+ x_2^3+x_3^3+x_4^3}$}} \author{C.~Baciu\thanks{The first
    and third author were supported by the DFG--Schwerpunkt ``Globale Methoden
    in der komplexen Geometrie''}\and V.~Ene\thanks{The second author is
    grateful to Oldenburg University for support.}\and G.~Pfister$^*$\and
  D.~Popescu\thanks{The fourth author was supported mainly by MSRI in Berkeley
    and by DFG in Kaiserslautern, partially he was also supported by Eager,
 and
    the Romanian Grants: EURROMMAT, CNCSIS and CERES Contracts 152/2001 and
39/2002. He is grateful to all these
    institutions.}}

\maketitle \thispagestyle{empty} \setcounter{page}{0}

\tableofcontents

\begin{abstract}We describe, by matrix factorizations, all the rank two
  maximal Cohen--Macaulay modules over singularities of type $x_1^3+
  x_2^3+x_3^3+x_4^3$.
\end{abstract}

\section*{Introduction}
Let $R$ be a hypersurface ring, that is $R=S/(f)$ for a regular
local ring $(S,\mathfrak{m})$ and $0\neq f\in \mathfrak{m}$.
Accordingly to  Eisenbud \cite{Ei}, any maximal Cohen--Macaulay
(briefly MCM) module over $R$ has a minimal free resolution of
periodicity $2$ which is completely given by a matrix
factorization $(\varphi,\psi), \varphi, \psi$ being square
matrices over $S$ such that $\varphi\psi=\psi\varphi=f\id_n$, for a
certain positive integer $n$. Therefore, in order to describe the
MCM $R$--modules, it is enough to describe their matrix
factorizations. In this paper we give the description, by matrix
factorizations, of the graded, rank two, indecomposable, MCM
modules over $K[x_1,x_2,x_3,x_4]/(x_1^3+x_2^3+x_3^3+x_4^3)$. Part
of this study was done with the help of the Computer Algebra
System \textsc{Singular} \cite{GPS}.\\

The MCM modules over the hypersurface $f_3=x_1^3+x_2^3+x_3^3$ were
described in \cite{LPP} as $1$--parameter families indexed by the
points of the curve $Z=V(f_3)\subset\mathbf{P}^2$. This
description is mainly based on the Atiyah's theory of the vector
bundles classification over elliptic curves, in particular over
$Z,$ and on difficult computations made with the Computer Algebra
System SINGULAR. The description depends on two discrete
invariants --- the rank and the degree of the bundle --- and on a
continuous invariant --- the points of the curve $Z$.\\
It is of high interest the classification of vector bundles, in
particular of ACM bundles (i.e. those which corresponds to MCM
modules) over the singularities of higher dimension. In the paper
\cite{EP} are described the matrix factorizations which define the
graded MCM modules of rank one over $f_4=x_1^3+x_2^3+x_3^3+x_4^3$.
There is a finite number of such modules which correspond to $27$
lines, $27$ pencils of  quadrics and $72$ nets of twisted cubic
curves lying on the surface $Y=V(f_4)\subset \mathbb{P}^3$. From
geometrical point of view the problem is easy, but the effective
description of the matrix factorizations is
difficult and SINGULAR has been intensively used.\\
In the present paper we continue this study for the graded MCM
modules of rank two. We obtain a general description of the MCM
orientable modules of rank two. They are given by skew--symmetric
matrix factorizations (see Theorem \ref{thm2.1}). The technique is
based on the results of Herzog and K\"uhl (see \cite{HK})
concerning the so called {\em Bourbaki exact sequences}. The
matrix factorizations of the graded, orientable, rank two,
$4$--generated MCM modules are parameter families indexed by the
points of the surface $Y$, that is two parameter families and some
finite ones in bijection with rank one MCM modules described in
\cite{EP} (see Theorem \ref{thm3.2}, here an important fact is
that  two Gorenstein ideals of codimension $2$ define the same MCM
module via the associated Bourbaki  sequence if and only if they
belong to  the same even linkage class). We also describe the
non--orientable MCM modules of rank two  over $f_4$. There is a
finite number of such modules which correspond somehow to the rank
one  modules described in \cite{EP}. The graded MCM modules,
non--orientable, of rank two are $2$--syzygy over $f_4$ of some
ideals of the form $J/(f_4),$ $J $ being an ideal of the
polynomial  ring $S=K[x_1,x_2,x_3,x_4],$ ($K$ is an algebraically
closed field of characteristic zero), with $f_4\in J,$ $\dim
S/J=2,\ \depth S/J=1$, whose Betti numbers over $S$ satisfy
$\beta_1(J)=\beta_0(J)+1$ and $\beta_2(J)=1$ (see Lemma
\ref{test}). This result has been essential in the description of
the graded, non--orientable MCM modules. The paper highlights
bijections between the classes of indecomposable, graded,
non--orientable MCM modules of rank two, $4$ and $5$--generated
and the classes of rank one, graded, MCM modules (see Theorem
\ref{prop4.6} and Theorem \ref{thm5.3}). Consequently, there
exists a bijection between the classes of indecomposable, graded,
non--orientable MCM modules of rank three,  $5$--generated and the
classes of rank one, graded, MCM modules (see Corollary
\ref{cor5.4}).
 These results remind us the theory of
Atiyah and give a small hope that the non-orientable case behaves
in the same way for higher rank. We also show that there are no
indecomposable, graded, non--orientable MCM modules of rank two
$6$--generated. Consequently, there exist no indecomposable,
graded, non--orientable MCM modules of rank four, 6-generated.

Till now the description of graded rank two MCM modules is not too
far from the theory of Atiyah. But the description of graded, rank
two, 6-generated MCM modules is different (see Section 6) as we
expected since a part of them given by Gorenstein ideals defined
by 5 general points on $Y$ forms a 5-parameter family (see
\cite{Mi}, \cite{IK}). However we believe that behind these
results there exists a nice theory of graded MCM modules over a
cubic hypersurface in four variables which waits to be discovered.

We express our thanks to A.~Conca, R.~Hartshorne, J.~Herzog and G.~Valla for
very helpful discussions on Section 6 and Theorem \ref{thm2.1}.

\section{Preliminaries}

Let $R_n:=K[x_1,x_2,\ldots,x_n]/(f_n)$, where
$f_n=x_1^3+x_2^3+\ldots+x_n^3$ and $K$ is an algebraic closed
field of characteristic $0$. Using the classification of vector
bundles over elliptic curves obtained by Atiyah \cite{At}, Laza,
Pfister and Popescu \cite{LPP} describe the matrix factorizations
of the graded, indecomposable and reflexive modules over $R_3$.
They give  canonical normal forms for the matrix factorizations of
all graded reflexive $R_3$--modules of rank one (see Section 3 in
\cite{LPP}) and show effectively how we can produce the
indecomposable graded reflexive $R_3$--modules of rank $\geq 2$
using \textsc{Singular} (see Section 5 in \cite{LPP}). We recall
from \cite{LPP} the description of the rank one, three-generated,
non--free, graded MCM $R_3$--modules since we shall use it in the
last section of our paper. First we recall the notations. Let
$P_0=[-1:0:1]\in V(f_3)$. For each
$\lambda=[\lambda_1:\lambda_2:1]\in V(f_3),\lambda\neq P_0$, we
set
$$\alpha_{\lambda}=\left(%
\begin{array}{ccc}
  0 & x_1-\lambda_1x_3 & x_2-\lambda_2x_3 \\
  x_1+x_3 & -x_2-\lambda_2x_3 & -wx_3 \\
  x_2 & wx_3 & (1-\lambda_1)x_3-x_1 \\
\end{array}%
\right),$$ where $w=\frac{\lambda_2^2}{\lambda_1+1}$ and, if
$\lambda=[\lambda_1:1:0]\in V(f_3)$, we set
$$\alpha_{\lambda}=\left(%
\begin{array}{ccc}
  0 & x_1-\lambda_1x_2 & x_3 \\
  x_1+x_3 & -\lambda_1x_1 & \lambda_1x_1+\lambda_1^2x_2 \\
  x_2 & x_3-x_1 & -x_1 \\
\end{array}%
\right).$$ Let $\beta_{\lambda}$ the adjoint matrix of
$\alpha_{\lambda}$.

\begin{thm} [(3.7) in \cite{LPP}]
$(\alpha_{\lambda},\beta_{\lambda})$ is a matrix factorization for
all $\lambda\in V(f_3), \lambda\neq P_0$, and the set of
three--generated MCM graded $R_3$--modules,
$$\mathcal{M}_0=\{\Co \alpha_{\lambda}\ |\ \lambda\in V(f_3), \lambda\neq
P_0\}$$ has the following properties:
\begin{enumerate}
    \item [(i)] All the modules from $\mathcal{M}_0$ have rank one.

    \item [(ii)] Every two different modules from $\mathcal{M}_0$ are not isomorphic
    \item [(iii)] Every three--generated, rank one, non--free,
    graded MCM $R_3$--module is isomorphic with one module from
    $\mathcal{M}_0$.
\end{enumerate}
\end{thm}

Now we consider the case $n=4$. In this case we do not have the
support of Atiyah classification. The complete description by
matrix factorizations of the rank one, graded, indecomposable MCM
modules over $R_4$ was given in \cite{EP}. \\
The aim of the present paper is to classify the rank two, graded,
indecomposable MCM modules over $R_4$. From now on, we shall
denote $R=R_4$, $f=f_4$ and we preserve the hypothesis on $K$ to
be algebraically closed and of characteristic zero.\\
Let $M$ be a rank two MCM module over $R$ and let $\mu(M)$ be the
minimal number of generators of $M$. By Corollary 1.3 of
\cite{HK}, we obtain that $\mu(M)\in \{3,4,5,6\}$. \\
First of all we consider {\bf the three--generated case}.
 The description of the rank one MCM $R$--modules is given in
\cite{EP}. We recall the notations. For $a,b,c,d,\varepsilon\in K$
such that $a^3=b^3=c^3=d^3=-1, \varepsilon^3=1,\varepsilon\neq 1$,
and $bcd=\varepsilon a$, we set
$$\alpha(b,c,d,\varepsilon)=\left(%
\begin{array}{ccc}
  0 & x_1-ax_4 & x_2-bx_3 \\
  x_1-cx_2 & -b^2x_3-abc^2\varepsilon^2x_4 & b^2c^2x_3-abc\varepsilon^2x_4 \\
  x_3-dx_4 & c^2x_2+bc^2x_3+acx_4 & -x_1-cx_2-ax_4 \\
\end{array}%
\right)$$ and
$$\beta(b,c,d,\varepsilon)=\alpha(b,c,d,\varepsilon)^t,$$
that is, the transpose of $\alpha(b,c,d,\varepsilon)$. Then each
of the matrices $\alpha(b,c,d,\varepsilon)$ and
$\beta(b,c,d,\varepsilon)$ forms with its adjoint,
$\alpha(b,c,d,\varepsilon)^{\ast}$, respectively
$\beta(b,c,d,\varepsilon)^{\ast}$, a matrix factorization of $f$. \\
For $a,b,c\in K$, distinct roots of $-1$, and $\varepsilon$ as
above, we set
$$\eta(a,b,c,\varepsilon)=\left(%
\begin{array}{ccc}
  0 & x_1+x_2 & x_3-ax_4 \\
  x_1+\varepsilon x_2 & -x_3+cx_4 & 0 \\
  x_3-bx_4 & 0 & -x_1-\varepsilon^2x_2 \\
\end{array}%
\right)$$ and
$$\vartheta(a,b,c)=\left(%
\begin{array}{ccc}
  0 & x_1+x_3 & x_2-ax_4 \\
  x_1-a^2bx_3 & -x_2+cx_4 & 0 \\
  x_2-bx_4 & 0 & -x_1+ab^2x_3 \\
\end{array}%
\right).$$ The matrices $\eta(a,b,c,\varepsilon)$ and
$\vartheta(a,b,c)$ form with their adjoint,
$\eta(a,b,c,\varepsilon)^{\ast}$, respectively
$\vartheta(a,b,c)^{\ast}$, a matrix factorization of $f$.

\begin{thm}[(3.4) in \cite{EP}] Let
$$\mathcal{M}=\{\Co
  \alpha(b,c,d,\varepsilon), \Co \beta(b,c,d,\varepsilon)\ |\
  b,c,d,\varepsilon \in K,$$
  $$
  b^3=c^3=d^3=-1,bcd=\varepsilon a, \varepsilon^3=1,\varepsilon\neq 1\}$$
  and
  $$\mathcal{N}=\{\Co \eta(a,b,c,\varepsilon), \Co\vartheta(a,b,c)\ |\
  \varepsilon^3=1,\varepsilon\neq 1 $$
  $$\mbox{\ and\ } (a,b,c) \mbox {\ is\
    a\ permutation\ of\ the\ roots\ of\ }-1 \}.$$
  Then the sets $\mathcal{M},
  \mathcal{N}$ of rank one, three--generated, MCM graded $R$--modules have the
  following properties:
 \begin{enumerate}
 \item [(i)] Every three--generated, rank one, indecomposable, graded MCM $R$--module
   is isomorphic with one module from $\mathcal{M}\cup \mathcal{N}$.
 \item [(ii)] If $M=\Co \alpha(b,c,d,\varepsilon)$ (or $M=\Co
   \beta(b,c,d,\varepsilon)$) belongs to $\mathcal{M}$ and $N\in \mathcal{M}$,
   then $N\simeq M$ if and only if $N=\Co
   \alpha(b\varepsilon,c\varepsilon,d\varepsilon,\varepsilon^2)$ (or $N=\Co
   \beta(b\varepsilon,c\varepsilon,d\varepsilon,\varepsilon^2)$).
    \item [(iii)] Any two different modules from $\mathcal{N}$ are not isomorphic.
    \item [(iv)] Any module of $\mathcal{N}$ is not isomorphic with some module of $\mathcal{M}$.
 \end{enumerate}
\end{thm}
The map $M\mapsto \Omega_R^1(M)$ is a bijection between the
three--generated, indecomposable, graded, MCM $R$--modules of rank
two and the three--generated, indecomposable, graded, MCM
$R$--modules of rank one. Thus, from the above theorem we obtain
the description of the rank two, three--generated, indecomposable,
graded MCM $R$--modules.

\begin{thm} Let
  $$\mathcal{M}^{\ast}=\{\Co \alpha(b,c,d,\varepsilon)^{\ast}, \Co
  \beta(b,c,d,\varepsilon)^{\ast}\ |\ b,c,d,\varepsilon \in K,$$
  $$
  b^3=c^3=d^3=-1,bcd=\varepsilon a, \varepsilon^3=1,\varepsilon\neq 1\}$$
  and
  $$\mathcal{N}^{\ast}=\{\Co \eta(a,b,c,\varepsilon)^{\ast},
  \Co\vartheta(a,b,c)^{\ast}\ |\ \varepsilon^3=1,\varepsilon\neq 1 $$
  $$\mbox{\ and\ } (a,b,c) \mbox {\ is\ a\ permutation\ of\ the\ roots\ of\
  }-1 \}.$$
  Then the sets $\mathcal{M}^{\ast}, \mathcal{N}^{\ast}$ of rank
  two, three--generated, MCM graded $R$--modules have the following properties:
 \begin{enumerate}
 \item [(i)] Every three--generated, rank two, indecomposable, graded MCM
   $R$--module is isomorphic with one module from $\mathcal{M}^{\ast}\cup
   \mathcal{N}^{\ast}$.

 \item [(ii)] If $M=\Co \alpha(b,c,d,\varepsilon)^{\ast}$ (or $M=\Co
   \beta(b,c,d,\varepsilon)^{\ast}$) belongs to $\mathcal{M}^{\ast}$ and $N\in
   \mathcal{M}^{\ast}$, then $N\simeq M$ if and only if $N=\Co
   \alpha(b\varepsilon,c\varepsilon,d\varepsilon,\varepsilon^2)^{\ast}$ (or
   $N=\Co
   \beta(b\varepsilon,c\varepsilon,d\varepsilon,\varepsilon^2)^{\ast}$).

 \item [(iii)] Any two different modules from $\mathcal{N}^{\ast}$ are not
   isomorphic.

\item [(iv)] Any module of $\mathcal{N}^{\ast}$ is not isomorphic
with some
  module of $\mathcal{M}^{\ast}$.
 \end{enumerate}
\end{thm}

\begin{cor}
There are $72$ isomorphism classes of rank two, indecomposable,
graded MCM modules over $R$ with three generators.
\end{cor}

\section{Skew symmetric matrices and rank 2 orientable MCM modules}

Let $\varphi=(a_{ij})_{1\leq i,j\leq 2s}$ be a generic skew
symmetric matrix, that is
$$a_{ii}=0, a_{ij}=-a_{ji}, \mbox{\ for\
all\ } i,j=\overline{1,2s}.$$ Then
$$\det(\varphi)=\Pf
(\varphi)^2,$$ where $\Pf(\varphi)$ denotes the Pfaffian of
$\varphi$ (see \cite[\S 5, no.\ 2]{Bo1} or \cite[(3.4)]{BH}). Like
determinants Pfaffians can be developed along a row. Set
$\varphi_{ij}$ the matrix obtained from $\varphi$ by deleting the
$i^{\mbox{th}}$ and $j^{\mbox{th}}$ rows and columns. Then, for
all $i=1,\ldots,2s,$

\begin{equation}\label{eq1}
\Pf (\varphi)=\sum_{\substack{j=1\\ j\neq
i}}^{2s}(-1)^{i+j}\sigma(i,j)a_{ij}\Pf(\varphi_{ij}),
\end{equation}

where $\sigma(i,j)$ denotes $\sgn (j-i)$. Multiplying (\ref{eq1})
by $\Pf(\varphi)$, we get

\begin{equation}\label{eq2}
\det(\varphi)=\sum_{j=1}^{2s}(-1)^{i+j}a_{ij}b_{ij},
\end{equation}

for $b_{ij}=\sigma(i,j)\Pf(\varphi_{ij})\Pf(\varphi)$ when $i\neq
j$ and $b_{ii}=0$. Since $\varphi$ is a generic matrix we see from
(\ref{eq2}) that $b_{ij}$ is exactly the algebraic complement of
$a_{ij}$ and so the transpose matrix $B$ of $(b_{ij})$ is the
adjoint matrix of $\varphi$. Set
$$\psi=\frac{1}{\Pf(\varphi)}B.$$
Then
$$\varphi\psi=\psi\varphi=\Pf(\varphi)\id_{2s},$$
as it is stated also in [\cite{JP}, \S 3].

\begin{prop}\label{prop2.1}
Let $f=x_1^3+x_2^3+x_3^3+x_4^3$ and $\varphi$ a skew symmetric
matrix over $S=K[x_1,x_2,x_3,x_4]$ of order $4$ or $6$ such that
$\det \varphi=f^2$, K being a field. Then $\Co \varphi$ is a MCM
module over $R:=S/(f)$ of rank $2$.
\end{prop}

\begin{proof}
Let $\psi$ be given for $\varphi$ as above, that is the $(i,j)$
entry of $\psi$ is $\sigma(i,j)\Pf(\varphi_{ij})$. As above we
have
$$\varphi\psi=\psi\varphi=f\cdot \id_n,\ n=4 \mbox{\ or\ }6$$
because $\Pf(\varphi)=f$. Then $(\varphi,\psi)$ is a matrix
factorization which defines a MCM $R$--module of rank $2$.
\end{proof}

\begin{thm}\label{thm2.1}
Preserving the hypothesis of Proposition \ref{prop2.1}, the
cokernel of a homogeneous skew symmetric matrix over $S$ of order
$4$ or  $6$ of determinant $f^2$ defines a graded  MCM $R$--module
$M$ of rank two. Conversely, each non--free graded orientable MCM
$R$--module $M$ of rank two is the cokernel of a map given by a
skew symmetric homogeneous matrix $\varphi$ over $S$ of order $4$
or  $6$, whose determinant is $f^2$ and $\varphi$ together with
$\psi$, defined above, form the matrix factorization of $M$.
\end{thm}

\begin{proof}
After Herzog and K\"{u}hl \cite{HK}, $M$ must be $4$ or $6$
minimally generated. Suppose that $M$ is $6$--generated (the other
case is similar). Then $M$ is the second syzygy over $R$ of a
Gorenstein ideal $I\subset R$ of codimension $2$ which is
$5$--generated by \cite{HK}. Using Buchsbaum--Eisenbud Theorem
(see e.g.\ \cite{BH}, (3,4)) there exists an exact sequence

\begin{equation}\label{eq3}
0\arrow{e}
S(-5)\arrow{e,t}{d_3}S^5(-3)\arrow{e,t}{d_2}S^5(-2)\arrow{e,t}{d_1}S
\end{equation}

such that $J=\im d_1, I=J/(f), d_2$ is a skew symmetric
homogeneous matrix, $d_3$ is the dual of $d_1$, $d_3=d_1^t$, and
$$d_1=\Bigl(\Pf \bigl((d_2)_1\bigr),-\Pf \bigl((d_2)_2\bigr),\ldots, \Pf
\bigl((d_2\bigr)_5)\Bigr),$$ where $(d_2)_i$ denotes the $4\times
4$ skew symmetric matrix obtained by deleting the $i^{\mbox{th}}$
row and column of $d_2$. Since $f\in J$ there exists $v:S(-1)\To
S^5 $ such that $d_1 v=f$ ($v$ is given by linear forms). It is
easy to see from (\ref{eq3}) that $I=J/(f)$ has the following
minimal resolution over $S:$
$$0\arrow{e}S(-5)\arrow{e,t}{\big(\substack{d_3\\0}\big)}S^6(-3)\arrow{e,t}{(d_2,v)}S^5(-2)
\arrow{e,t}{\tilde{d}_1}I\arrow{e}0.$$ Like in \cite{Ei}, since
$fI=0$, there exists a map $h:S^5(-5)\ra S^6(-3)$ such that $(d_2,
v)h=f\cdot \id_5$ and we get the following exact sequence

\begin{equation}\label{eq4}
R^6(-5)\arrow{e,t}{\big(\bar{h}|\substack{\bar{d}_3\\
0}\big)}R^6(-3)\arrow{e,t}{(\bar{d}_2,\bar{v})}R^5(-2)\arrow{e,t}{\bar{d}_1}
I\arrow{e}0.
\end{equation}
On the other hand, $\varphi=\left(%
\begin{array}{cc}
  d_2 & v \\
  -v^t & 0 \\
\end{array}%
\right)$ is a skew symmetric homogeneous matrix of order $6$. Let
$\psi$ given as above. By construction $\psi$ has the form $\left(%
\begin{array}{cc}
  C & -d_1^t \\
  d_1 & 0 \\
\end{array}%
\right)$ and so $(d_2,v)\big(\substack{C\\ d_1}\big)=f\cdot
\id_5$. Taking $h=\big(\substack{C\\ d_1}\big)$ above, we get from
(\ref{eq4}), the following exact sequence:

$$R^6(-6)\arrow{e,t}{\varphi}R^6(-5)\arrow{e,t}{\psi}
R^6(-3)\arrow{e,t}{(\bar{d}_2,\bar{v})}
R^5(-2)\arrow{e,t}{\bar{d}_1}I\arrow{e}0,$$

which gives
$$\Co \varphi\cong \im \psi =\Omega_R^2(I).$$
We have $\varphi\big(h|\substack{-d_1^t\\0}\big)=f\cdot \id_6$ and
so $\det(\varphi)$ is a power of $f$. Since the entries of
$\varphi$ are linear forms, we get $\det(\varphi)=f^2$.
\end{proof}

\section{Orientable, rank 2, 4--generated MCM modules}

Let $K$ be an algebraically closed field of characteristic zero,\\
$S=K[x_1,x_2,x_3,x_4], f=x_1^3+x_2^3+x_3^3+x_4^3$, and $R=S/(f)$.
Let $M$ be a graded, indecomposable, $4$--generated MCM
$R$--module of rank two. After Herzog and K\"{u}hl \cite{HK},
$M\cong \Omega_R^2(I)$, where $I$ is a graded $3$--generated
Gorenstein ideal such that $\dim R/I=1$. Then $I=J/(f)$, with
$J\subset S$ a graded, $3$--generated ideal containing $f$. Let
$\alpha_1,\alpha_2,\alpha_3$ be a minimal system of homogeneous
generators of $J$. Since $\dim S/J=1$, it follows  that
$\alpha_1,\alpha_2,\alpha_3$ is a regular system of elements in
$S$.\\

Let $u,a,b\in K$ with $a^3=b^3=-1, u^2+u+1=0$ and $\sigma=(i\ j\
s)$ be a permutation of the set $\{2,3,4\}$ with $i<j$. Set
$$w_{\sigma 1}=x_1-ax_s,\ w_{\sigma 2}=x_i-bx_j,$$ $$v_{\sigma
1}=x_1^2+ax_1x_s+a^2x_s,\ v_{\sigma 2}=x_i^2+bx_ix_j+b^2x_j^2.$$
then we have
$$f=w_{\sigma 1}v_{\sigma 1}+w_{\sigma 2}v_{\sigma 2}.$$ Let
$\lambda=[\lambda_1:\lambda_2:\lambda_3:1]$ be a point of the
surface $V(f)\subset \mathbb{P}^3$. We set
$$p_{i\lambda}=x_i-\lambda_ix_4, \mbox{\ and\ }
q_{i \lambda}=x_i^2+\lambda_i x_ix_4+\lambda_i^2 x_4^2,\mbox{\
for\ }1\leq i\leq 3.$$ Let $\lambda=[\lambda_1:\lambda_2:1:0]$ be
a point of $V(f)$. We set $$p_{i\lambda}=x_i-\lambda_ix_3, q_{i
\lambda}=x_i^2+\lambda_i x_ix_3+\lambda_i^2 x_3^2,\mbox{\ for\
}1\leq i\leq 2$$ and $$p_{3\lambda}=x_4,\ q_{3\lambda}=x_4^2.$$ If
$\lambda=[\lambda_1:1:0:0]\in V(f)$, we set
$$p_{1\lambda}=x_1-\lambda_1x_2,\ q_{1 \lambda}=x_1^2+\lambda_1
x_1x_2+\lambda_1^2 x_2^2 $$ and $$ p_{2\lambda}=x_3,\
p_{3\lambda}=x_4, q_{2\lambda}=x_3^2,\ q_{3\lambda}=x_4^2.$$ In
all cases we have
$$f=\sum_{i=1}^3p_{i\lambda}q_{i\lambda}.$$
These are the only ways to write
$f$ as a linear combination of two or three forms of degrees $1, 2$, provided
that the 1--forms are linearly independent.  Since $f\in
(\alpha_1,\alpha_2,\alpha_3)$, we may suppose that either $\alpha_i$ is in the
set $\{p_{i\lambda},q_{i\lambda}\}$ for each $1\leq i\leq 3$, or $\alpha_i$ is
in the set $\{w_{\sigma i},v_{\sigma i}\}$ for each $1\leq i\leq 2$ and
$\beta=\alpha_3$ is a regular element in $R/(\alpha_1,\alpha_2)$.

\begin{lem}\label{lem3.1}
Let $M$ be a graded, indecomposable, $4$--generated MCM
$R$--module of rank $2$. Then $M$ is one of the following modules:
\begin{itemize}
    \item [(1)] $\Omega_R^2(p_{1 \lambda},p_{2 \lambda},p_{3
    \lambda})$ or $\Omega_R^2(q_{1 \lambda},q_{2 \lambda},q_{3
    \lambda})$, for some $\lambda\in V(f),$
    \item [(2)] $\Omega_R^2(w_{\sigma 1},v_{\sigma 2},\beta)$
 or  $\Omega_R^2(w_{\sigma 2},v_{\sigma 1},\beta)$ for some
 $a,b,\sigma$and $\beta$ as above.
\end{itemize}
\end{lem}

\begin{proof}
Set $$I_{\lambda}= (p_{1 \lambda},p_{2 \lambda},p_{3 \lambda})$$
and $$\varphi_{\lambda}=\left(%
\begin{array}{cccc}
  0 & p_{3\lambda} & -p_{2\lambda} & -q_{1\lambda} \\
  -p_{3\lambda} & 0 & -p_{1\lambda} & q_{2\lambda} \\
  p_{2\lambda} & p_{1\lambda} & 0 & q_{3\lambda} \\
  q_{1\lambda} & -q_{2\lambda} & -q_{3\lambda} & 0 \\
\end{array}%
\right),\ \psi_{\lambda}=\left(%
\begin{array}{cccc}
  0 & -q_{3\lambda} & q_{2\lambda} & p_{1\lambda} \\
  q_{3\lambda} & 0 & q_{1\lambda} & -p_{2\lambda} \\
  -q_{2\lambda} & -q_{1\lambda} & 0 & -p_{3\lambda} \\
  -p_{1\lambda} & p_{2\lambda} & p_{3\lambda} & 0 \\
\end{array}%
\right).$$ We have the following exact sequence:
$$R^3(-5)\oplus
R(-6)\arrow{e,t}{\varphi_{\lambda}}R^4(-4)\arrow{e,t}{\psi_{\lambda}}
R^3(-2)\oplus R(-3)\arrow{e,t}{A}R^3(-1)\arrow{e,t}{\tau}
I_{\lambda}\arrow{e}0,$$ where
$\tau=(-p_{1\lambda},p_{2\lambda},p_{3\lambda})$ and $A$ is given
by the first three rows of $\varphi_\lambda$. Thus
$\Omega^2(I_{\lambda})\cong \Co (\varphi_{\lambda})$ and
$(\varphi_{\lambda},\psi_{\lambda})$ is a matrix factorization of
$\Omega^2(I_{\lambda})$. The ideals $I_{\lambda}$ and $(q_{1
\lambda},q_{2 \lambda}, p_{3 \lambda})$ belong to the same even
linkage class since
$$I_{\lambda}\sim (q_{1 \lambda},p_{2 \lambda},
p_{3 \lambda})\sim (q_{1 \lambda},q_{2 \lambda}, p_{3 \lambda}).$$
For the first link we consider the regular sequence $\{p_{1
\lambda}q_{1 \lambda}, p_{2 \lambda}, p_{3 \lambda}\}$ and for the
second one the sequence $\{q_{1 \lambda}, p_{2 \lambda}q_{2
\lambda}, p_{3 \lambda}\}$. Similarly one can see that
$I_{\lambda}$ is evenly linked with the ideals $(q_{1 \lambda},
p_{2 \lambda}, q_{3 \lambda})$ and $(p_{1 \lambda}, q_{2 \lambda},
q_{3 \lambda})$. By Theorem 2.1 \cite{HK}, we get that
$$\Co(\varphi_{\lambda})\cong\Omega_R^2(I_{\lambda})\cong\Omega_R^2(q_{1\lambda
},q_{2\lambda},p_{3\lambda}) \cong \Omega_R^2(q_{1\lambda
},p_{2\lambda},q_{3\lambda})\cong \Omega_R^2(p_{1\lambda
},q_{2\lambda},q_{3\lambda}).$$ Analogously we see that $$\Co
(\psi_{\lambda})\cong\Omega_R^2(q_{1\lambda
},q_{2\lambda},q_{3\lambda})\cong \Omega_R^2(p_{1\lambda
},p_{2\lambda},q_{3\lambda})\cong\Omega_R^2(p_{1\lambda
},q_{2\lambda},p_{3\lambda})\cong $$ $$
\cong\Omega_R^2(q_{1\lambda },p_{2\lambda},p_{3\lambda}).$$ Thus
the case when $\alpha_i$ is one of the forms
$\{p_{i\lambda},q_{i\lambda}\}$ gives $(1)$.\\
Now let $\sigma, a,b$ as above and $\beta\in S$ which is regular
on $R/(w_{\sigma 1},v_{\sigma 2})$. Set $$I_{\sigma
\beta}(a,b,u)=(w_{\sigma 1},v_{\sigma 2},\beta)$$ and $$\varphi_{\sigma \beta}(a,b,u)=\left(%
\begin{array}{cccc}
  0 & w_{\sigma 1} & -v_{\sigma 2} & 0 \\
  -w_{\sigma 1} & 0 & -\beta & w_{\sigma 2} \\
  v_{\sigma 2} & \beta & 0 & v_{\sigma 1} \\
  0 & -w_{ \sigma2} & -v_{\sigma 1} & 0 \\
\end{array}%
\right),$$ $$ \psi_{\sigma \beta}(a,b,u)=\left(%
\begin{array}{cccc}
  0 & -v_{\sigma 1} & w_{\sigma 2} & \beta \\
  v_{\sigma 1} & 0 & 0 & -v_{\sigma 2} \\
  -w_{\sigma 2} & 0 & 0 & -w_{\sigma 1} \\
  -\beta & v_{\sigma 2} & w_{\sigma 1} & 0 \\
\end{array}%
\right).$$ We have the following exact sequence:
$$\begin{diagram}
\node{R^4}\arrow{e,t}{\varphi_{\sigma \beta}(a,b,u)}\node{R^4}
\arrow{e,t}{\psi_{\sigma \beta}(a,b,u)}
\node{R^4}\arrow{e,t}{B}\node{R^3}\arrow{e,t}{\tau'}\node{I_{\sigma
\beta}(a,b,u)}\arrow{e}\node{0}
\end{diagram},$$
where $\tau'=(-\beta,v_{\sigma 2},w_{\sigma 1})$ and $B$ is the
matrix given by the first three rows of $\varphi_{\sigma
\beta}(a,b,u)$. Thus $$\Omega_R^2\bigl(I_{\sigma
\beta}(a,b,u)\bigr)\cong \Co \bigl(\varphi_{\sigma
\beta}(a,b,u)\bigr).$$ As above we see that
$$\Omega_R^2\bigl(I_{\sigma \beta}(a,b,u)\bigr)\cong \Omega_R^2(w_{\sigma
2},v_{\sigma 1},\beta)$$ and
$$\Omega_R^2(w_{\sigma 1},w_{\sigma 2},\beta)\cong
\Omega_R^2(v_{\sigma 1},v_{\sigma 2},\beta)\cong \Co
\bigl(\psi_{\sigma \beta}(a,b,u)\bigr).$$ Thus the case when
$\alpha_i$ is one of the forms $\{w_{\sigma i},v_{\sigma i}\}$ for
$i\leq 2$ gives  $(2)$.
\end{proof}

Let $$\mathcal{M}=\{\Co(\varphi_{\lambda}),\Co (\psi_{\lambda})\
|\ \lambda\in V(f)\}.$$ For $a,b,\sigma$ as above, set
$$\varphi_{\sigma}(a,b,u)=\varphi_{\sigma,x_jx_s}(a,b,u),\
\psi_{\sigma}(a,b,u)=\psi_{\sigma,x_jx_s}(a,b,u),$$ that is
$\beta=x_jx_s$. Let
$$\mathcal{P}=\{\Co\bigl(\varphi_{\sigma}(a,b,u)\bigr),\Co(\psi_{\sigma}(a,b,u))\
|\ a,b,\sigma \mbox{\ as\ above\ }\}.$$

\begin{thm}\label{thm3.2} The set
$\mathcal{M}\cup\mathcal{\mathcal{P}}$contains only
non--isomorphic, indecomposable, graded, orientable,
$4$--generated MCM $R$--modules of rank 2 and every
indecomposable, graded, orientable, $4$--generated MCM $R$--module
of rank 2 is isomorphic with one module of
$\mathcal{M}\cup\mathcal{\mathcal{P}}$.
\end{thm}

\begin{proof}
  Applying Lemma \ref{lem3.1} we must show in the case $(2)$ that $\beta$ can
  be taken $x_jx_s$. Since $v_{\sigma 1}-w_{\sigma 1}(x_1+2ax_s)=3a^2x_s^2,$
  adding in $\varphi_{\sigma\beta}(a,b,u)$ multiples of the last row to the
  second one and multiples of the first column to the third one, we may
  suppose the entry $(2,3)$ of the form $\gamma+x_s\delta$, with
  $\gamma,\delta$ depending only on $x_j,x_i$. These transformations modify
  the entries $(2,2),(3,3)$ which are now possibly non--zero. Adding similar
  multiples of the last column to the second one and multiples of the first
  row to the third one, we get $\varphi_{\sigma,\beta}(a,b,u)$ of the same type
  as before but with $\beta=\gamma+x_s\delta$. We may reduce to consider $\delta\not \in K.$
  Indeed, if $\delta\in K,$ then, acting on the rows and columns of $\varphi_{\sigma\beta}(a,b,u),$
  we get that $M=\Co(\varphi_{\sigma\beta}(a,b,u))$ is decomposable or belongs to the set
  $\mathcal{M}.$ Now let $\delta$ be not constant. Similarly, adding in
  $\varphi_{\sigma\beta}(a,b,u)$ multiples of the first row to the second one
  and multiples of the last column to the third one we may suppose that the
  entry $(2,3)$ has the form $\varepsilon x_jx_s$ with $\varepsilon\in K$.
  These transformations modify the entries $(2,2),(3,3)$. After similar
  transformations we get $\varphi_{\sigma\beta}(a,b,u)$ of the same type as
  before but with $\beta=\varepsilon x_jx_s$. If $\varepsilon=0$ we see that
  $\varphi_{\sigma\beta}(a,b,u)$ is a direct sum of two $2\times 2$--matrices
  which contradicts the indecomposability of
  $M=\Co\bigl(\varphi_{\sigma\beta}(a,b,u)\bigr)$. So $\varepsilon\neq 0$.
  Divide the second and the third column of $\varphi_{\sigma\beta}(a,b,u)$ with
  $\varepsilon$ and multiply the first and the last row of
  $\varphi_{\sigma\beta}(a,b,u)$ with $\varepsilon$. We reduce to the
  case $\varepsilon=1$, that is $\beta=x_jx_s$.\\
  Now we show that two different modules from \mbox{$\mathcal{M}\cup
    \mathcal{P}$} are not isomorphic. Note that the Fitting ideals of
  \mbox{$\varphi_{\lambda}$} (respectively \mbox{$\psi_{\lambda}$}) modulo
  \mbox{$(x_1,\ldots,x_4)^2$} have the form \mbox{$(p_{1 \lambda},p_{2
      \lambda},p_{3 \lambda})$} and the Fitting ideals of
  \mbox{$\varphi_{\sigma}(a,b,u)$} (respectively \mbox{$\psi_{\sigma}(a,b,u)$})
  modulo \mbox{$(x_1,\ldots,x_4)^2$} have the form \mbox{$(w_{\sigma
      1},w_{\sigma 2})$} and these ideals are all different. Thus
$$\{\Co(\varphi_{\lambda})\ |\ \lambda\in V(f)\}\cup \{\Co
    \bigl(\varphi_{\sigma}(a,b,u)\bigr)\ |\ \sigma, a, b \mbox{\ as\  above\ }\}$$ contains
  only non--isomorphic modules (similarly for $\psi$--es). It follows that, if
  $N,P\in \mathcal{M}\cup \mathcal{P}$ are isomorphic and different, then
  $N\simeq \Omega_R^1(P)$. \\ If $N=\Co (\varphi_{\lambda})$, for $\lambda\in
  V(f)$, then this is not possible since the ideals $(p_{1 \lambda}, p_{2
    \lambda},p_{3 \lambda})$ and $(q_{1 \lambda}, q_{2 \lambda},q_{3
    \lambda})$ are not in the same even linkage class. Indeed, by the proof of
  $(i)$ in Lemma \ref{lem3.1}, $(p_{1 \lambda}, p_{2 \lambda},p_{3 \lambda})$
  is evenly linked with $(q_{1 \lambda}, q_{2 \lambda},p_{3 \lambda})$ and
  this last ideal is obviously directly linked with $(q_{1 \lambda}, q_{2
    \lambda},q_{3 \lambda})$. If $N=\Co\bigl(\varphi_{\sigma}(a,b,u)\bigr)$ for
  some $\sigma, a,b$, and $N\simeq \Omega_R^1(N)$, then the ideals $(w_{\sigma
    1},v_{\sigma 2}, x_jx_s)$ and $(w_{\sigma 1},w_{\sigma 2}, x_jx_s)$ are
  evenly linked. But these ideals are directly linked by the regular sequence
  $\{w_{\sigma 1},v_{\sigma 2}w_{\sigma 2}, x_jx_s\}$, contradiction!\\ It
  remains to show that $\mathcal{M}\cup \mathcal{P}$ contains only
  indecomposable modules. If $N\in \mathcal{M}$, let us say
  $N=\Co(\varphi_{\lambda})$ for $\lambda=[\lambda_1:\lambda_2:\lambda_3:1]$
  we see that $N/x_4N$ is exactly the module corresponding to the matrix $$
\left(%
\begin{array}{cccc}
  0 & x_3 & -x_2 & -x_1^2 \\
  -x_3 & 0 & -x_1 & x_2^2 \\
  x_2 & x_1 & 0 & x_3^2 \\
  x_1^2 & -x_2^2 & -x_3^2 & 0 \\
\end{array}%
\right)$$ whose cokernel is the  special module $M_2$ (see
\cite{LPP} for the special module of rank two which corresponds to
the special bundle from Atiyah classification). Thus $N/x_4N$ is
indecomposable and, by Nakayama's Lemma, $N$ is indecomposable.
Now let $N\in \mathcal{P},$
$N=\Co\bigl(\psi_{\sigma}(a,b,u)\bigr)$. By the permutation of the
rows and the columns of $\psi_{\sigma}(a,b,u)$, we may suppose
that it has the form:
$$\left(%
\begin{array}{cccc}
  w_{\sigma 1} & -v_{\sigma 2} & x_jx_s & 0 \\
  w_{\sigma 2} & v_{\sigma 1} & 0 & x_jx_s \\
  0 & 0 & v_{\sigma 1} & v_{\sigma 2} \\
  0 & 0 & -w_{\sigma 2} & w_{\sigma 1} \\
\end{array}%
\right).$$ Suppose $N$ is decomposable. Then
$\psi_{\sigma}(a,b,u)$ is a direct sum of two matrices of order
two which we may suppose to be given by the submatrices of the
above one given by the first two lines and columns respectively
the last two lines
 and columns (this is obvious modulo $x_j$ or $x_s$).
  Due to the particular form of $\psi_{\sigma}(a,b,u)$ this
means that there exist two matrices $A,B$ of order two such that
$$x_jx_s\cdot \id_2=\left(\begin{array}{cc}
  w_{\sigma 1} & -v_{\sigma 2} \\
  w_{\sigma 2}& v_{\sigma 1} \\
\end{array}\right)A+B\left(%
\begin{array}{cc}
  v_{\sigma 1} & v_{\sigma 2} \\
   -w_{\sigma 2} & w_{\sigma 1} \\
\end{array}%
\right)$$ which is impossible.
\end{proof}

\begin{rems}\label{rems3.3}
(i) There exists a bijection between $$\mathcal{P}_1=\{\Co
\bigl(\varphi_{\sigma}(a,b,u)\bigr)\ |\ \sigma, a,b\}$$ and the
two--generated non--free MCM $R$--modules which remind us Atiyah's
classification. Thus $\mathcal{P}_1$ contains $54$ modules
corresponding to $27$ lines and $27$ pencils of conics of $V(f)$.
Similarly, $\mathcal{P}_2=\{\Co
\bigl(\psi_{\sigma}(a,b,u)\bigr)\ |\ \sigma, a,b\}$ contains $54$ modules.\\
(ii) $\mathcal{M}$ is a kind of ``blowing up'' of
$M_2,\Omega_R^1(M_2)$ from \cite{LPP} (see the proof of Theorem
\ref{thm3.2}). Note also that $\mathcal{M}$ consists of two
classes of modules parameterized by the points of
$V(f)$, which is also in Atiyah's idea.\\
(iii) The matrices $\varphi$ defining the modules of
$\mathcal{M}\cup \mathcal{P}$ are skew symmetric as our Theorem
\ref{thm2.1} predicted.
\end{rems}

\section{Non--orientable, rank 2, 4--generated MCM modules}

Let $M$ be a graded non--orientable, rank 2, MCM $R$--module,
without free direct summands. We should like to express $M$ as a
2--syzygy of an ideal $I$, $M\cong \Omega_R^2(I)$, with
$\mu(M)=\mu(I)+1$ (this is known in orientable case by \cite{HK},
see here Section 3).

The following proposition can be found in \cite[Korollar 2]{B}.

\begin{prop}\label{prop4.3} Let $(A,m)$ be a Noetherian normal local
domain with $\dim A\geq 2$ and $N$ a finite torsion--free
$A$--module. Then there exists a finite free submodule $F\subset
N$ such that $N/F$ is isomorphic with an ideal of $A$ and the
canonical  map $F/mF\rightarrow N/mN$ is injective.
\end{prop}


Applying Proposition \ref{prop4.3} we obtain the following exact sequence:
\begin{equation}\label{sir1}
 0\rightarrow R\rightarrow M \rightarrow I \rightarrow 0
\end{equation}
for an ideal $I\subset R$, which induces an exact sequence

 $$0\rightarrow K=R/m\rightarrow M/mM \rightarrow I/mI
\rightarrow 0.$$

Thus $\mu(M)=\mu(I)+1$. \\

As we know in the orientable case to get MCM $R$--modules of rank
2 we must choose $I$ such that Ext$^1_R(I,R)$ is a cyclic
$R$--module or more precisely such that $R/I$ is Gorenstein. In
the non--orientable case one can also show that Ext$^1_R(I,R)$
must be a cyclic $R$--module, but this is not very helpful since
it is hard to check this condition for arbitrary $I$. Below we
shall state an easier condition.

Let $J\subset S=K[X_1,\ldots,X_4]$ be an ideal such that $f\in J$
and $I=J/(f)$.

\begin{lem}\label{test}
Let

$$0\arrow{e} S^{s_3}\arrow{e,t}{d_3} S^{s_2}\arrow{e,t}{d_2} S^{s_1}
\arrow{e,t}{d_1} J \arrow{e} 0$$

be a minimal free $S$--resolution of an ideal $J$ with depth
$S/J=1$. \\ If $\Rk \Omega^2_R\bigl(J/(f)\bigr)=2$ and
$\mu(\Omega^2_R\bigl(J/(f)\bigr)=\mu(I)+1$ then $s_1=s_2\leq 5$
and $s_3=1$.

\end{lem}

\begin{proof} As in the proof of Theorem \ref{thm2.1} we get a minimal free resolution of
$I=J/(f)$ over $S$ in the following way:

Let $v:S\arrow{e} S^{s_1}$ be an $S$--linear map such that
$jd_1v=f\id_S$, where $j:J\arrow{e} S$ is the inclusion. Let
${\tilde d}_1 $ be the composite map $S^{s_1}\arrow{e,t}{d_1}
J\arrow{e} J/(f)=I$. Then the following sequence

$$0\arrow{e} S^{s_3}\arrow{e,t}{\big(\substack{d_3\\0}\big)} S^{s_2+1}
 \arrow{e,t}{(d_2,v)} S^{s_1}\arrow{e,t}{\tilde {d}_1} I\arrow{e} 0$$

is exact and forms a minimal free $S$--resolution of $I$ over $S$.
Since $$f\cdot S^{s_1}\subset \im(d_2,v),$$ there exists an
$S$--linear map $h:S^{s_1} \arrow{e} S^{s_2+1}$ such that
$$(d_2,v)h=f\id_{S^{s_1}} $$ and we get the following exact sequence

$$R^{s_3+s_2}\arrow{e,t}{({\bar h}|\substack{{\bar d}_3\\0})} R^{s_2+1}
\arrow{e,t}{({\bar d}_2,{\bar v})} R^{s_1}\arrow{e,t}{{\bar d}_1}
I \arrow{e} 0,$$

which is part of a minimal free $R$--resolution of $I$. Thus
$M=\Omega_R^2(I)$ is the image of the first map above and so
$s_2+s_3=s_2+1=s_1+1$ because
$\mu(M)=\mu\bigl(\Omega_R^1(M)\bigr)=\mu(I)+1$ by hypothesis. It
follows $s_3=1$, $s_1=s_2$. As $\mu(M)\leq 3 \Rk_RM=6$
 we get $s_1\leq 5$.
\end{proof}

Let $\det N$ be the corresponding class of the bidual
$(\wedge^{n}N)^{**},\ n=\Rk N$, in $Cl(R)$ for a torsion free
$R$--module $N$. Since $\det$ is an additive function, we get
$\det(M)=0$ if and only if $\det(I)=0$. Thus $M$ is
non--orientable if and only if $I$ is non--orientable, that is
$\codim (J)\leq 1$ for all ideals $J\subset R$ isomorphic
with $I$,
after \cite{HK}. Since $M$ has rank $2$, we get
$\codim (I)=1$. Thus $\dim R/I=2$ and, from (\ref{sir1}), we get
$\depth R/I=1$, that is $R/I$ is not Cohen--Macaulay. Also from
(\ref{sir1}) we get $\Omega_R^2(M)\simeq \Omega_R^2(I)$ and so
$M\simeq \Omega_R^2(I)$.

\begin{prop}\label{prop4.4}
Each graded, non--orientable, rank two, $s$--generated MCM
$R$--module is the second syzygy $\Omega_R^2(I)$ of an
$(s-1)$--generated graded ideal $I\subset R$ with $\depth R/I=1$
and $\dim R/I=2$.
\end{prop}

As in Section $3$, let $u,a,b\in K$, with $$a^3=b^3=-1,\
u^2+u+1=0,$$ $\sigma=(i\ j\ s)$ be a permutation of the set
$\{2,3,4\}$ with $i<j$ and set $$w_{\sigma 1}=x_1-ax_s, w_{\sigma
2}=x_i-bx_j,$$ $$ v_{\sigma 1}=x_1^2+ax_1x_s+a^2x_s^2, v_{\sigma
2}=x_i^2+bx_ix_j+b^2x_j^2.$$ We have $$ v_{\sigma 1}= v'_{\sigma
1} v''_{\sigma 1}, \ v_{\sigma 2}= v'_{\sigma 2} v''_{\sigma 2}$$
for $$v'_{\sigma 1}=x_1-uax_s,\  v''_{\sigma 1}=x_1+(1+u)ax_s,$$
$$v'_{\sigma 2}=x_i-ubx_j,\  v''_{\sigma 2}=x_i+(1+u)bx_j.$$ Set
$$I_{1\sigma}(a,b,u)=(x_sv'_{\sigma 2},v_{\sigma 2},w_{\sigma 1}),$$
$$I_{2\sigma}(a,b,u)=(x_jv''_{\sigma 1},v_{\sigma 1},w_{\sigma 2}),$$
$$I_{3\sigma}(a,b,u)=(x_sv''_{\sigma 2},v_{\sigma 2},v_{\sigma 1}),$$
$$I_{4\sigma}(a,b,u)=(x_jv'_{\sigma 1},v_{\sigma 1},v_{\sigma 2}).$$

Set $$\varphi_{1\sigma}(a,b,u)=\left(%
\begin{array}{cccc}
  0 & w_{\sigma 1} & -v''_{\sigma 2} & 0 \\
  -w_{\sigma 1} & 0 & -x_s & w_{\sigma 2} \\
  v_{\sigma 2} & x_sv'_{\sigma 2} & 0 & v_{\sigma 1} \\
  0 & -w_{\sigma 2}v'_{\sigma 2} & -v_{\sigma 1} & 0 \\
\end{array}%
\right),$$

$$\psi_{1\sigma}(a,b,u)=\left(%
\begin{array}{cccc}
  0 & -v_{\sigma 1} & w_{\sigma 2} & x_s \\
  v_{\sigma 1} & 0 & 0& -v''_{\sigma 2} \\
  -w_{\sigma 2}v'_{\sigma 2} & 0 & 0 & -w_{\sigma 1} \\
 -x_sv'_{\sigma 2} & v_{\sigma 2} & w_{\sigma 1} & 0 \\
\end{array}%
\right),$$

$$\varphi_{2\sigma}(a,b,u)=\left(%
\begin{array}{cccc}
  0 & w_{\sigma 2} & -v'_{\sigma 1} & 0 \\
  -w_{\sigma 2} & 0 & -x_j & w_{\sigma 1} \\
  v_{\sigma 1} & x_jv''_{\sigma 1} & 0 & v_{\sigma 2} \\
  0 & -w_{\sigma 1}v''_{\sigma 1} & -v_{\sigma 2} & 0 \\
\end{array}%
\right),$$

$$\psi_{2\sigma}(a,b,u)=\left(%
\begin{array}{cccc}
  0 & -v_{\sigma 2} & w_{\sigma 1} & x_j \\
  v_{\sigma 2} & 0 & 0& -v'_{\sigma 1} \\
  -w_{\sigma 1}v''_{\sigma 1} & 0 & 0 & -w_{\sigma 2} \\
 -x_jv''_{\sigma 1} & v_{\sigma 1} & w_{\sigma 2} & 0 \\
\end{array}%
\right),$$

$$\varphi_{3\sigma}(a,b,u)=\left(%
\begin{array}{cccc}
  0 & v_{\sigma 1} & -v'_{\sigma 2} & 0 \\
  -v_{\sigma 1} & 0 & -x_s & w_{\sigma 2} \\
  v_{\sigma 2} & x_sv''_{\sigma 2} & 0 & w_{\sigma 1} \\
  0 & -w_{\sigma 2}v''_{\sigma 2} & -w_{\sigma 1} & 0 \\
\end{array}%
\right),$$

$$\psi_{3\sigma}(a,b,u)=\left(%
\begin{array}{cccc}
  0 & -w_{\sigma 1} & w_{\sigma 2} & x_s \\
  w_{\sigma 1} & 0 & 0& -v'_{\sigma 2} \\
  -w_{\sigma 2}v''_{\sigma 2} & 0 & 0 & -v_{\sigma 1} \\
 -x_sv''_{\sigma 2} & v_{\sigma 2} & v_{\sigma 1} & 0 \\
\end{array}%
\right),$$

$$\varphi_{4\sigma}(a,b,u)=\left(%
\begin{array}{cccc}
  0 & v_{\sigma 2} & -v''_{\sigma 1} & 0 \\
  -v_{\sigma 2} & 0 & -x_j & w_{\sigma 1} \\
  v_{\sigma 1} & x_jv'_{\sigma 1} & 0 & w_{\sigma 2} \\
  0 & -w_{\sigma 1}v'_{\sigma 1} & -w_{\sigma 2} & 0 \\
\end{array}%
\right),$$

$$\psi_{4\sigma}(a,b,u)=\left(%
\begin{array}{cccc}
  0 & -w_{\sigma 2} & w_{\sigma 1} & x_j \\
  w_{\sigma 2} & 0 & 0& -v''_{\sigma 1} \\
  -w_{\sigma 1}v'_{\sigma 1} & 0 & 0 & -v_{\sigma 2} \\
 -x_jv'_{\sigma 1} & v_{\sigma 1} & v_{\sigma 2} & 0 \\
\end{array}%
\right).$$

\begin{thm}\label{prop4.6}
\begin{itemize}
\item [(i)] For each $1\leq t\leq 4$, the pair
  $\bigl(\varphi_{t \sigma}(a,b,u),\psi_{t \sigma}(a,b,u)\bigr)$ forms a matrix
  factorization of $\Omega_R^2\bigl(I_{t \sigma}(a,b,u)\bigr)$.
    \item [(ii)] The set

      $$\mathcal{N}=\{\Co\bigl(\varphi_{t \sigma}(a,b,u)\bigr),
      \Co\bigl(\psi_{t \sigma}(a,b,u)\bigr)\ |\ 1\leq t\leq 4,\ \sigma, a,
      b, u
      \}$$

contains only graded, indecomposable,
    non--orientable, $4$--generated MCM $R$--modules of rank $2$.
    \item [(iii)] Every indecomposable, graded, non--orientable,
$4$--generated MCM
    module over $R$ of rank $2$ is isomorphic with one module of
$\mathcal{N}$.
\item [(iv)] The modules of $\mathcal{N}$ are pairwise non--isomorphic. In particular,
there exist $432$ isomorphism classes of  indecomposable, graded, non--orientable,
$4$--generated MCM module over $R$ of rank $2.$
\end{itemize}
\end{thm}

\begin{proof}
$(i)$. It is easy to check that
$$\varphi_{t \sigma}(a,b,u)\cdot\psi_{t \sigma}(a,b,u)=f\cdot \id_4$$
and the following sequence is exact:
$$\begin{diagram}
\node{R(-6)^4}\arrow{e,t}{\varphi_{1\sigma}(a,b,u)}\node{R(-5)^2\oplus
R(-4)^2}\arrow{e,t}{\psi_{1\sigma}(a,b,u)}\node{R(-3)^4}\arrow{e,t}{A_1}\node{}\\
\node{}\arrow{e,t}{A_1}\node{R(-2)^2\oplus
R(-1)}\arrow{e}\node{I_{1\sigma}(a,b,u)}\arrow{e}\node{0,}
\end{diagram}$$
where $A_1$ is the $3\times 4-$matrix formed by the first three
rows of $\varphi_{1 \sigma}(a,b,u)$. Thus $(i)$ holds for $t=1$,
the other
cases being similar.\\
$(ii)$. Clearly $I_{1\sigma}(a,b,u)\subset (v'_{\sigma
2},w_{\sigma 1})$ and so $\dim R/I_{1\sigma}(a,b,u)=2$. As $x_s$
is zero--divisor in $R/I_{1\sigma}(a,b,u)$ we see that $\depth
R/I_{1\sigma}(a,b,u)=1$ and, by Proposition \ref{prop4.5},
$\Omega_R^2(I)$ is non--orientable, $4$--generated of rank $2$.
Note that after some linear transformations
$\varphi_{1\sigma}(a,b,u)$ becomes

$$\left(%
\begin{array}{cccc}
  w_{\sigma 2} &  -w_{\sigma 1} & 0 & x_s \\
  v_{\sigma 1} & v_{\sigma 2}   & x_s v'_{\sigma 2} & 0 \\
  0 & 0 & w_{\sigma 1} & v''_{\sigma 2} \\
  0 & 0 & -w_{\sigma 2}v'_{\sigma 2} & v_{\sigma 1} \\
\end{array}%
\right),$$

and as in the last part of the proof of Theorem \ref{thm3.2} we
see that $\Co\bigl(\varphi_{1\sigma}(a,b,u)\bigr)$ is
indecomposable because there exist no two matrices $A,B$ of order
two such that

$$\left(\begin{array}{cc}
  0 & x_s \\
  x_s v'_{\sigma 2} & 0 \\
 \end{array}\right)
=\left(\begin{array}{cc}
  w_{\sigma 2} & -w_{\sigma 1} \\
   v_{\sigma 1} &  v_{\sigma 2}\\
\end{array}\right)A+B\left(%
\begin{array}{cc}
  w_{\sigma 1} & v''_{\sigma 2} \\
   -w_{\sigma 2}v'_{\sigma 2} & v_{\sigma 1} \\
\end{array}%
\right).$$

Similarly follows the cases $t>1$.\\
$(iii)$. Now let $M$ be an indecomposable, graded,
non--orientable, $4$--generated MCM $R$--module of rank $2$. By
Proposition \ref{prop4.4} there exists a graded ideal $I\subset R$
with $\dim R/I=2,\ \depth R/I=1$, which is $3$--generated and such
that $M\simeq \Omega_R^2(I)$. Then $I=J/(f)$ with $J\subset
S=K[x_1,x_2,x_3,x_4]$ a three generated ideal containing $f$. Let
$\alpha_1,\alpha_2,\alpha_3$ be a minimal system of homogeneous
generators of $J$. If $f$ does not belong to the ideal generated
by two $\alpha_t$, then, as in Section $3$, $f=\sum_{t=1}^3p_tq_t$
and, after a renumbering, we may suppose that $\alpha_t$ is
necessarily either $p_t$ or $q_t$, for all $1\leq t\leq 3$. Then
$\alpha_1,\alpha_2,\alpha_3$ is a regular system of elements in
$S$ and so $R/I=S/I$ is Cohen--Macaulay which is false.\\
Thus we may suppose $f\in (\alpha_1,\alpha_2)$. Then there exist
$a,b\in K$ with $a^3=b^3=-1$ and $\sigma=(i\ j\ s)$ a permutation
of the set $\sigma=\{2,3,4\},\ i<j$, such that $\alpha_t$ is
necessarily either $w_{\sigma t}$ or $v_{\sigma t}$, for $t=1,2$.
If $\alpha_1=w_{\sigma 1}, \alpha_2=w_{\sigma 2}$, then
$R/(\alpha_1,\alpha_2)$ is a domain and
$\alpha_1,\alpha_2,\alpha_3$ must be a regular system of elements
in $S$ and so, again, $R/I=S/I$ is Cohen--Macaulay,
contradiction!\\
We have the following cases:
\begin{center}
{\bf Case I:\;  $\alpha_1=w_{\sigma 1}$ }
\end{center}
Then $\alpha_2$ must be $v_{\sigma 2}$ and we have
$$(\alpha_1,\alpha_2)=(v'_{\sigma 2},w_{\sigma 1})\cap (v''_{\sigma
2},w_{\sigma 1}).$$ It follows that a zero--divisor of
$R/(\alpha_1,\alpha_2)$ must be either in $(v'_{\sigma
2},w_{\sigma 1})$ or in $(v''_{\sigma 2},w_{\sigma 1})$. As we
know $\alpha_3$ is a zero--divisor in $R/(\alpha_1,\alpha_2)$ and
so $\alpha_3\in (v'_{\sigma 2},w_{\sigma 1})$ or $\alpha_3\in
(v''_{\sigma
2},w_{\sigma 1})$. \\
{\bf I (a).} Suppose  $$\alpha_3\in (v'_{\sigma 2},w_{\sigma
1}).$$ Subtracting from $\alpha_3$ a multiple of $w_{\sigma 1},$
we may take $\alpha_3=v'_{\sigma 2}\beta$ for a form $\beta$ of
$S$. Note that the matrices
$$\varphi=\left(%
\begin{array}{cccc}
  0 & w_{\sigma 1} & -v''_{\sigma 2} & 0 \\
  -w_{\sigma 1} & 0 & -\beta & w_{\sigma 2} \\
  v_{\sigma 2} & \beta v'_{\sigma 2} & 0 & v_{\sigma 1} \\
  0 & -w_{\sigma 2}v'_{\sigma 2} & -v_{\sigma 1} & 0 \\
\end{array}%
\right),$$

$$\psi=\left(%
\begin{array}{cccc}
  0 & -v_{\sigma 1} & w_{\sigma 2} & \beta \\
  v_{\sigma 1} & 0 & 0& -v''_{\sigma 2} \\
  -w_{\sigma 2}v'_{\sigma 2} & 0 & 0 & -w_{\sigma 1} \\
 -\beta v'_{\sigma 2} & v_{\sigma 2} & w_{\sigma 1} & 0 \\
\end{array}%
\right)$$  give the following exact sequence:
$$\arrow{e}R^4\arrow{e,t}{\varphi}R^4\arrow{e,t}{\psi}R^4\arrow{e,t}{B_1}R^3\arrow{e}I
\arrow{e} 0,$$ where $B_1$ is given by the first three rows of
$\varphi$. Thus $(\varphi,\psi)$ is a matrix factorization of
$\Omega_R^2(I)\simeq M$. Adding in $\varphi$ multiples of the
first row to the second one and adding multiples of the forth
column to the third one, we may suppose that the entry $(2,3)$ of
$\varphi$ depends only on $x_1,x_s$. These transformations modify
also the entries $(2,2)$ and $(3,3)$ which are now not zero.
Adding similar multiples of first column to the second one and of
the forth row to the third one, we get $\varphi$ of the same type
as before but with $\beta$ depending only on $x_1,x_s$. Since
$v_{\sigma 1}-w_{\sigma 1}(x_1+2ax_s)=3ax_s^2$, adding in
$\varphi$ multiples of the first column to the third one and
multiples of the forth row to the second row, we may suppose that
the entry $(2,3)$ has the form $\lambda x_s$ for some $\lambda \in
K$. These transformations modify also the entries $(3,3)$ and
$(2,2)$ which are now not zero. Adding similar multiples of the
first row to the third one and of the fourth column to the second
column, we get $\varphi$ of the same type as before but with
$\beta=\lambda x_s$. If $\lambda =0$, then clearly $\varphi$ is
the direct sum of two $2$--matrices which contradicts that $M$ is
indecomposable. So $\lambda\neq 0$. Now we divide the second and
the third column of $\varphi$ by $\lambda$ and multiply the first
and the fourth row by $\lambda$. The new $\varphi$ is as before
but with $\lambda=1$, that is $\varphi=\varphi_{1 \sigma}(a,b,u)$.\\
{\bf I (b).} Suppose $$\alpha_3\in (v''_{\sigma 2},w_{\sigma
1}).$$ Then we may take $\alpha_3=v''_{\sigma 2} \beta$, for a
 form $\beta$. With  a similar proof as above, we obtain $M\simeq
\Co \bigl(\psi_{3 \sigma}(a,b,u)\bigr)$.

\begin{center}
{\bf Case II:\; $\alpha_2=w_{\sigma 2}$.}
\end{center}
 Then $\alpha_1=v_{\sigma 1}$. It follows that
$(\alpha_1,\alpha_2)=(v'_{\sigma 1},w_{\sigma 2})\cap (v''_{\sigma
1},w_{\sigma 2})$. We have the following two
subcases: \\
{\bf II (a).} $\alpha_3\in (v'_{\sigma 1},w_{\sigma 2})$. We may
suppose $\alpha_3=v'_{\sigma 1}\beta$, for a form $\beta$
and we obtain that $M\simeq \Co\bigl(\psi_{4 \sigma}(a,b,u)\bigr)$.\\
{\bf II (b).} $\alpha_3\in (v''_{\sigma 1},w_{\sigma 2})$. In this
subcase we may  take $\alpha_3=v''_{\sigma 1}\beta$, for  a form
$\beta$ and we obtain that $M\simeq \Co\bigl(\varphi_{2
\sigma}(a,b,u)\bigr)$.

\begin{center}
{\bf Case III:\; $\alpha_1=v_{\sigma 1},\alpha_2=v_{\sigma 2} $.}
\end{center}
Then $(\alpha_1,\alpha_2)=(v'_{\sigma 1},v'_{\sigma 2})\cap
(v'_{\sigma 1},v''_{\sigma 2})\cap (v''_{\sigma 1},v'_{\sigma
2})\cap (v''_{\sigma 1},v''_{\sigma 2})$. We proceed like in the
above cases taking $\alpha_3$ from one prime ideal of the above
decomposition of $\bigl(\alpha_1,\alpha_2)\bigr)$, let us say
$\alpha_3\in (v'_{\sigma 1},v'_{\sigma 2}) $, that is
$\alpha_3=v'_{\sigma 1}\beta+v'_{\sigma 2}\gamma$ for some
$\beta,\gamma\in S$. Suppose that one cannot reduce the problem to
the case $\beta=0$ or $\gamma=0$, this implies for example that
$v'_{\sigma 1}$ does not divide $\gamma$ and $v'_{\sigma 2}$ does
not divide $\beta$. Then
$\Omega^1_S\bigl((\alpha_1,\alpha_2,\alpha_3)\bigr) \subset S^3$
contains the columns of the following matrix

$$\left(%
\begin{array}{cccc}
v_{\sigma 2} & \alpha_3 & 0 & v''_{\sigma 2}\beta \\
- v_{\sigma 1} & 0 & \alpha_3 & v''_{\sigma 1}\gamma \\
0  & - v_{\sigma 1} & - v_{\sigma 2} & -v''_{\sigma 1}v''_{\sigma 2}\\
\end{array}%
\right)$$

and we can see that
$\mu(\Omega^1_S\bigl((\alpha_1,\alpha_2,\alpha_3)\bigr)\geq 4$,
which contradicts Lemma \ref{test}. Thus we may suppose, let us
say $\alpha_3=v'_{\sigma 1}\beta$ where $\beta$ is not a multiple
of $v''_{\sigma 1}$. Now we may proceed as in the above cases and
we obtain, in order, $M\simeq \Co\bigl(\varphi_{4
\sigma}(a,b,u)\bigr), M\simeq \Co\bigl(\varphi_{3
\sigma}(a,b,u)\bigr), M\simeq \Co\bigl(\psi_{1
  \sigma}(a,b,u)\bigr)$, and $M\simeq \Co\bigl(\psi_{2 \sigma}(a,b,u)\bigr)$.\\
$(iv)$.  We shall prove that the matrices of the set
$$\mathcal{N}'=\{ \varphi_{t\sigma}(a,b,u),\
\psi_{t\sigma}(a,b,u)\ |\ 1\leq t\leq 4, \sigma, a, b, u\}$$ are
pairwise non--equivalent. We shall consider the matrices which are
obtained from the matrices of $\mathcal{N}'$ reducing their
entries modulo $\mathfrak{m}^2$.  If $A,B\in \mathcal{N}'$ are
equivalent, then there exist $P,Q$, two invertible $4\times
4$--matrices with the entries in $K[x_1,x_2,x_3,x_4]$ such that
$PA=BQ$.  Let $\widetilde{A}$ and $\widetilde{B}$ be the matrices obtained
from $A,$ respectively $B,$ by reducing modulo $\mathfrak{m}^2$ their
entries. From the equality $PA=BQ,$ we obtain that there exist two
invertible scalar matrices $\widetilde{P},\widetilde{Q}\in
\mathcal{M}_4(K)$ such that
$\widetilde{P}\widetilde{A}=\widetilde{B}\widetilde{Q}$.  This means that the
matrices $\widetilde{A},\widetilde{B}$ are also  equivalent by some scalar
invertible matrices. We construct the "reduced" matrices
$\tilde{\varphi}_{t\sigma}(a,b,u)$ and
$\tilde{\psi}_{t\sigma}(a,b,u),$ for all $t$.  We see that the
matrices
$\tilde{\varphi}_{1\sigma}(a,b,u),\tilde{\varphi}_{2\sigma}(a,b,u),
\tilde{\psi}_{3\sigma}(a,b,u)$ and $\tilde{\psi}_{4\sigma}(a,b,u)$
have the entries of the last two rows zero and the rest of the
matrices have the entries of the first two columns zero. First we
choose two matrices $\widetilde{A},\widetilde{B}$, one of them with the
last two rows zero and the other with the first two  columns zero.
Suppose that $\widetilde{A}\sim \widetilde{B}$.  It results that there  are
two invertible scalar $4\times 4$--matrices $U,V$ such that
$$\widetilde{A}U=V\widetilde{B}.$$ From this equality we get that  the
last two rows in the matrix $V\widetilde{B}$ are zero. Looking at the
four possibilities to choose the matrix  $\widetilde{B},$ we see that the
non--zero elements of the columns $3$ and $4$ in $\widetilde{B}$ are linear
independent.
Therefore the  last two rows in $V$ must be zero, contradicting
$V$ invertible.\\
Hence we  could find two  equivalent matrices in
the set $\mathcal{N}'$ only if both have the last two rows zero or
the first two columns zero. It is clear that we may reduce the
study of the equivalent matrices $\widetilde{A},\widetilde{B}$ which have
the last two rows zero. Let $U,V\in \mathcal{M}_{4 \times 4}(K)$
be invertible matrices such that $\widetilde{A}U=V\widetilde{B}$.
Let $$
\widetilde{A}=\left(%
\begin{array}{cc}
  A_1 & A_2 \\
  0 & 0 \\
\end{array}%
\right), \widetilde{B}=\left(%
\begin{array}{cc}
  B_1 & B_2 \\
  0 & 0 \\
\end{array}%
\right),         U=\left(%
\begin{array}{cc}
  U_1 & U_2 \\
  U_3 & U_4 \\
\end{array}%
\right),     V=\left(%
\begin{array}{cc}
  V_1 & V_2 \\
  V_3 & V_4 \\
\end{array}%
\right),    $$ be the decomposition of our matrices in $2 \times
2$ blocks. We may suppose that $A_1$, $B_1$ are of the form $w_{\sigma i}\id_2$, $1\leq i\leq 2$.
Then $$A_2U_3=V_1B_1-A_1U_1.$$ If $A_1$ has on the main
diagonal the element $w_{\sigma 1},$ then $A_2$ has on the main
diagonal two elements from the set $\{w_{\sigma 2},-v'_{\sigma 2},
-v''_{\sigma 2}\}$.  Inspecting the elements in the above equality,
we get that, if $B_1$ has $w_{\sigma 2}$ on the main diagonal,
then $U_1=0$ and so $U_3$ is invertible. Then $A_2=V_1BU_3^{-1}=w_{\sigma 2}
V_1U_3^{-1}$ which is not possible.
 This means that it remains to study the
cases $$\widetilde{A}=\tilde{\varphi}_{1\sigma}(a,b,u),
\widetilde{B}=\tilde{\psi}_{3\tau}(n,p,v)$$ and
$$\widetilde{A}=\tilde{\varphi}_{1\sigma}(a,b,u),\widetilde{B}=\tilde{\varphi}_{1\tau}(n,p,v),$$
for some $\sigma, a,b,u,\tau,n,p,v$.  Let $U,V\in \mathcal{M}_{4
\times 4}(K)$ be  invertible  matrices such that
$$\tilde{\varphi}_{1\sigma}(a,b,u)\cdot U=V\cdot
\tilde{\psi}_{3\tau}(n,p,v).
$$  Comparing the elements of the first  row in
the above equality, we obtain that  $U$ has all the entries of the
third row zero, contradicting $U$ invertible.\\ In the same way we
check that if $\tilde{\varphi}_{1\sigma}(a,b,u)$ and
$\tilde{\varphi}_{1\tau}(n,p,v)$ are different, then they are not
equivalent.
\end{proof}

Let $M(\sigma, a,b) = \Co\bigl(\varphi_{\sigma}(a,b)\bigr)$,
$N(\sigma,a,b) = \Co\bigl(\psi_{\sigma}(a,b)\bigr)$ be the graded
rank one 2--generated MCM $R$--modules (see \cite{EP}).

\begin{rem}
  There exists an indecomposable extension in Ext$^1_R \bigl(M(\sigma,a,b),\;
  N(\tau,n,p)\bigr)$ if and only if $\sigma = \tau$.   In this case, there
  exists a unique indecomposable rank 2, 4--generated MCM module corresponding
  to the extension (up to an iso) which is orientable if $n = a, p = b$ and
  non--orientable otherwise.   Since all $N(\sigma,n,p)$ are 9, the result is
  that for fixed $M(\sigma,a,b)$ there exists just one orientable and eight
 non--orientable MCM--modules, which are extensions $E$ of the form
\[
0 \to N(\sigma,n,p) \to E \to M(\sigma,a,b)\to 0\,.
\]
So we have 27 orientable and $8 \times 27$ non--orientable MCM--modules.
Similarly, taking now extensions $F$ of the form
\[
0 \to M(\sigma,n,p) \to F \to N(\sigma,a,b)\to 0
\]
we obtain another 27 orientable and $8 \times 27$ non--orientable MCM--modules.
\end{rem}


\section{Non--orientable, rank 2, 5--generated MCM modules}

As in Section $3$, let $u,a,b\in K$, with $$a^3=b^3=-1,\
u^2+u+1=0,$$ $\sigma=(i\ j\ s)$ be a permutation of the set
$\{2,3,4\}$ with $i<j$ and set $$w_{\sigma 1}=x_1-ax_s, w_{\sigma
2}=x_i-bx_j,$$ $$ v_{\sigma 1}=x_1^2+ax_1x_s+a^2x_s^2, v_{\sigma
2}=x_i^2+bx_ix_j+b^2x_j^2.$$ We have $$ v_{\sigma 1}= v'_{\sigma
1} v''_{\sigma 1}, \ v_{\sigma 2}= v'_{\sigma 2} v''_{\sigma 2}$$
for $$v'_{\sigma 1}=x_1-uax_s,\  v''_{\sigma 1}=x_1+(1+u)ax_s,$$
$$v'_{\sigma 2}=x_i-ubx_j,\  v''_{\sigma 2}=x_i+(1+u)bx_j.$$

Consider the following ideals:

Set
$$J_{1\sigma}(a,b,u)=(v_{\sigma 1},v_{\sigma 2},v'_{\sigma 1}v'_{\sigma
2}, v''_{\sigma 1}v''_{\sigma 2}),$$
$$J_{2\sigma}(a,b,u)=(v_{\sigma 1},v_{\sigma 2},v'_{\sigma
1}v''_{\sigma 2}, v''_{\sigma 1}v'_{\sigma 2}).$$

Denote by ${\cal J}$ the union of the above families of ideals.

Set
$$\rho_{1\sigma}(a,b,u)=\left(%
\begin{array}{ccccc}
  0 & -v''_{\sigma 2} & -v'_{\sigma 2} &w_{\sigma 1} & 0\\
v'_{\sigma 1} & 0 & 0 &w_{\sigma 2} & -v''_{\sigma 2}v''_{\sigma 1}\\
- v''_{\sigma 2} & 0 & v''_{\sigma 1} & 0 & 0\\
0 & v'_{\sigma 1} & 0 & 0 &  v_{\sigma 2}\\
0 & -w_{\sigma 2}& 0 & 0 & w_{\sigma 1}v''_{\sigma 1}\\
\end{array}%
\right)$$

and

$$\omega_{1\sigma}(a,b,u)=\left(%
\begin{array}{ccccc}
 -w_{\sigma 2}v''_{\sigma 1} & w_{\sigma 1}v''_{\sigma 1} &
-w_{\sigma 2}v'_{\sigma 2} & 0 & v''_{\sigma 1}v''_{\sigma 2 }\\
0 & 0 & 0 & w_{\sigma 1}v''_{\sigma 1} & -v_{\sigma 2}\\
-w_{\sigma 2}v''_{\sigma 2} & -w_{\sigma 1}v''_{\sigma 2} &
w_{\sigma 1}
v'_{\sigma 1} & 0 & v_{\sigma 2 }^{\prime\prime 2}\\
v_{\sigma 1} & v_{\sigma 2} & v'_{\sigma 1}v'_{\sigma 2 } &
v''_{\sigma 1}v''_{\sigma 2 } & 0\\
0 & 0 & 0 & w_{\sigma 2} & v'_{\sigma 1}\\
\end{array}%
\right).$$

Clearly the pair of the matrices above forms the matrix
factorization of $\Omega^2_R\bigl(J_{1\sigma}(a,b,u)/(f)\bigr)$.
By permutations of $v'_{\sigma 1}, v''_{\sigma 1},v'_{\sigma 2},
v''_{\sigma 2}$ one can find easily the matrix factorization of
the module $\Omega^2_R\bigl(J_{2\sigma}(a,b,u)/(f)\bigr)$.

 Also set
$$T_{1\sigma}(a,b,u)=(v_{\sigma 1},v_{\sigma 2},
v'_{\sigma 1}v''_{\sigma 2},{v''_{\sigma 2}}^2),$$
$$T_{2\sigma}(a,b,u)=(v_{\sigma 1},v_{\sigma 2},
v''_{\sigma 1}v''_{\sigma 2},{v''_{\sigma 2}}^2),$$
$$T_{3\sigma}(a,b,u)=(v_{\sigma 1},v_{\sigma 2},
v''_{\sigma 1}v'_{\sigma 2},{v'_{\sigma 2}}^2),$$
$$T_{4\sigma}(a,b,u)=(v_{\sigma 1},v_{\sigma 2},
v'_{\sigma 1}v'_{\sigma 2},{v'_{\sigma 2}}^2),$$
$$T_{5\sigma}(a,b,u)=(v_{\sigma 1},v_{\sigma 2},
v'_{\sigma 1}v''_{\sigma 2},{v'_{\sigma 1}}^2),$$
$$T_{6\sigma}(a,b,u)=(v_{\sigma 1},v_{\sigma 2},
v'_{\sigma 1}v'_{\sigma 2},{v'_{\sigma 1}}^2),$$
$$T_{7\sigma}(a,b,u)=(v_{\sigma 1},v_{\sigma 2},
v''_{\sigma 1}v''_{\sigma 2},{v''_{\sigma 1}}^2),$$
$$T_{8\sigma}(a,b,u)=(v_{\sigma 1},v_{\sigma 2},
v''_{\sigma 1}v'_{\sigma 2},{v''_{\sigma 1}}^2),$$

 and denote by
${\cal T}$ the set of all these ideals.

Set $$\mu_{1 \sigma}(a,b,u)=\left(%
\begin{array}{ccccc}
  v''_{\sigma 2} & 0 & 0 & 0 & w_{\sigma 1} \\
  0 &  v''_{\sigma 2} & - v'_{\sigma 1} & 0 & w_{\sigma 2} \\
  - v''_{\sigma 1} & 0 &  v'_{\sigma 2} &  v''_{\sigma 2} & 0 \\
  0 & - v'_{\sigma 2} & 0 & - v'_{\sigma 1} & 0 \\
 0 & - v''_{\sigma 1}w_{\sigma 1} & 0 & w_{\sigma 2}v''_{\sigma 2} & 0 \\
\end{array}%
\right)$$ and $$\nu_{1 \sigma}(a,b,u)=\left(%
\begin{array}{ccccc}
   v'_{\sigma 2}w_{\sigma 2} & - v'_{\sigma 2}w_{\sigma 1}& - v'_{\sigma 1}w_{\sigma 1} & - v''_{\sigma 2}w_{\sigma 1} & 0 \\
  0& 0 & 0 & - v''_{\sigma 2}w_{\sigma 2} & -v'_{\sigma 1}\\
   v''_{\sigma 1}w_{\sigma 2} & - v''_{\sigma 1}w_{\sigma 1} & w_{\sigma 2}v''_{\sigma 2} &0 & -v''_{\sigma 2} \\
   0 & 0 & 0 & -v''_{\sigma 1}w_{\sigma 1}& v'_{\sigma 2}\\
v_{\sigma 1} &  v_{\sigma 2} &  v''_{\sigma 2} v'_{\sigma 1} & { v''_{\sigma 2}}^2 & 0 \\
\end{array}%
\right).$$

The pair of the matrices above forms the matrix factorization of
$\Omega^2_R\bigl(T_{1\sigma}(a,b,u)/(f)\bigr)$. By permutations of
$v'_{\sigma 1}, v''_{\sigma 1},v'_{\sigma 2}, v''_{\sigma 2}$ one
can find easily the matrix factorization for the $2$--syzygy of
the other ideals of ${\cal T}$.

\begin{lem}\label{lem5.1} Let $M$ be a graded non--orientable, rank two, 5--generated
 MCM $R$--module,
without free direct summands. Then there exists an ideal
 $J\in {\cal J}\cup {\cal T}$ such that $f\in J$ and  $M\cong
 \Omega^2\bigl(J/(f)\bigr)$. Conversely, for every $J\in {\cal J}\cup {\cal
   T}$,
the module $ \Omega^2\bigl(J/(f)\bigr)$ is a non--orientable, rank
two, 5--generated MCM $R$--module without free direct summands.
\end{lem}

\begin{proof} The second statement follows easily since we have
already the matrix factorizations above of those ideals. Let $M$
be as above. As in the beginning of Section 4 we see that $M\cong
\Omega^2\bigl(J/(f)\bigr)$, for $J$ an ideal of $S$ containing
$f$, with $\mu(J)=4$, dim $S/J=2$, depth $S/J$=1 and
$\mu\bigl(\Omega^1_S(J)\bigr)=5$. We may also suppose
$J=(\alpha_1, \alpha_2,\alpha_3,\alpha_4)$ with $f\in (\alpha_1,
\alpha_2)$, where $\alpha_t$ is necessarily either $w_{\sigma t}$
or $v_{\sigma t}$ for $t=1,2$ for some $a,b$ and a certain
permutation $\sigma $ as above. Clearly we
 cannot have simultaneously  $\alpha_t=w_{\sigma t}$ because then
$(\alpha_1, \alpha_2)$ is a prime ideal and one cannot find
$\alpha_3, \alpha_4$ zero divisors as we need. We treat the
following cases:

\begin{center}
{\bf Case I:}\; $\alpha_1=w_{\sigma 1}$
\end{center}

Then we have $\alpha_2=v_{\sigma 2}$ and $(\alpha_1, \alpha_2)$ is
the intersection of the prime ideals $(v'_{\sigma 2},w_{\sigma
1})$, $(v''_{\sigma 2},w_{\sigma 1})$. Since $\alpha_3, \alpha_4$
must be zero divisors in $S/(\alpha_3, \alpha_4)$ we have the
following possibilities:

(I1)  $\alpha_3=v'_{\sigma 2}\beta$, $ \alpha_4=v'_{\sigma
2}\gamma$, (I2)  $\alpha_3=v''_{\sigma 2}\beta$, $
\alpha_4=v''_{\sigma 2}\gamma$, (I3)  $\alpha_3=v'_{\sigma
2}\beta$, $ \alpha_4=v''_{\sigma 2}\gamma$, (I4)
$\alpha_3=v''_{\sigma 2}\beta$, $ \alpha_4=v'_{\sigma 2}\gamma$

for some homogeneous $\beta,\gamma$ from
$\mathfrak{m}=(x_1,x_2,x_3,x_4)$. In the first case we see that  the
relations given by the columns of the following matrix:

$$\left(%
\begin{array}{ccccc}
v_{\sigma 2}& \alpha_3 &\alpha_4   &0 & 0    \\
 -w_{\sigma 1} &0 & 0 & \gamma & \beta   \\
  0 &  -w_{\sigma 1} & 0 & 0 & -v''_{\sigma 2}\\
 0 & 0 & -w_{\sigma 1}  & -v''_{\sigma 2} & 0  \\
\end{array}%
\right),$$

are elements  in $\Omega^1_S(J)\subset S^4$. Clearly these columns
are part in the minimal system of generators of $\Omega^1_S(J)$
because $w_{\sigma 1},v''_{\sigma 2}$ form a regular system in
$S$. The subcase (I2) is similar, this contradicts Lemma
\ref{test}.

Suppose now (I3) holds. Then the relations given by the columns of
the following matrix

$$\left(%
\begin{array}{ccccc}

v_{\sigma 2}& v'_{\sigma 2}\beta    &v''_{\sigma 2}\beta & 0 & 0    \\
 -w_{\sigma 1} &0 & 0 & \beta & \gamma   \\
  0 &  -w_{\sigma 1} & 0 &  -v''_{\sigma 2} & 0\\
 0 & 0 & -w_{\sigma 1}  & 0 & -v'_{\sigma 2}  \\

\end{array}%
\right),$$

are part of a minimal set of generators of $\Omega^1_S(J)$ (note
that $w_{\sigma 1},v''_{\sigma 2},v'_{\sigma 2}$ form a regular
system in $S$). Contradiction! Case (I4) is similar.

\begin{center}
{\bf Case II:}\; $ \alpha_1=v_{\sigma 1}$, $ \alpha_2=v_{\sigma 2}$
\end{center}

Since $ (\alpha_1,\alpha_2)=(v'_{\sigma 1},v'_{\sigma 2})\cap
(v'_{\sigma 1},v''_{\sigma 2})\cap (v''_{\sigma 1},v''_{\sigma
2})\cap (v''_{\sigma 1},v'_{\sigma 2})$ we see that the zero
divisors of $S/ (\alpha_1,\alpha_2)$ must be in one of the prime
ideals of the above decomposition. Suppose $\alpha_3\in
(v'_{\sigma 1},v'_{\sigma 2})$. If $\alpha_3=\beta_1 v'_{\sigma
1}+ \beta_2 v'_{\sigma 2}$ then as in the proof of Case III of
Proposition \ref{prop4.6} we see that there are at least 4 minimal
 relations between first three $\alpha$. Then all $\alpha$ have at
least 5 minimal relations. Contradiction! Thus  $\alpha_3$ as well
$\alpha_4$ are multiples of one $v'_{\sigma t},v''_{\sigma t}$. So
we have the following
 possibilities:

(II1)   $\alpha_3=v'_{\sigma 1}\beta$, $ \alpha_4=v'_{\sigma
1}\gamma$, (II2)  $\alpha_3=v''_{\sigma 1}\beta$, $
\alpha_4=v''_{\sigma 1}\gamma$, (II3)  $\alpha_3=v'_{\sigma
2}\beta$, $ \alpha_4=v'_{\sigma 2}\gamma$, (II4)
$\alpha_3=v''_{\sigma 2}\beta$, $ \alpha_4=v''_{\sigma 2}\gamma$
(II5)  $\alpha_3=v'_{\sigma 1}\beta$, $ \alpha_4=v''_{\sigma
1}\gamma$ (II6)  $\alpha_3=v'_{\sigma 1}\beta$, $
\alpha_4=v''_{\sigma 2}\gamma$ (II7) $\alpha_3=v'_{\sigma 2}\beta,
\alpha_4=v''_{\sigma 1}\gamma$ (II8) $\alpha_3=v'_{\sigma 2}\beta,
\alpha_4=v''_{\sigma 2}\gamma$ (II9) $\alpha_3=v'_{\sigma 1}\beta,
\alpha_4=v'_{\sigma 2}\gamma$ (II10)$\alpha_3=v''_{\sigma 1}\beta,
\alpha_4=v''_{\sigma 2}\gamma$.

\begin{center}{\bf Subcase:}
 $\alpha_3=v'_{\sigma 1} \beta$, $\alpha_4=v'_{\sigma
1} \gamma$, $(v_{\sigma 2}v''_{\sigma 1},\gamma)\cong 1$,
$(v_{\sigma 2}v''_{\sigma 1}, \beta)\cong 1$
\end{center}

We see that the relations given by the columns of the following
matrix

$$\left(%
\begin{array}{ccccc}

v_{\sigma 2}& \beta    &\gamma & 0 & 0    \\
 -v_{\sigma 1} &0 & 0 & \alpha_3 & \alpha_4   \\
  0 &  -v''_{\sigma 1} & 0 &  -v_{\sigma 2} & 0\\
 0 & 0 & -v''_{\sigma 1}  & 0 & -v_{\sigma 2}  \\

\end{array}%
\right),$$

are part from a minimal system of generators of
 $\Omega^1_S(J)$ which must be false. Indeed, it is easy to see that
the last 4 columns are part in a minimal system of generators of
 $\Omega^1_S(J)$. If the first column belongs to the module generated
by the last four then there exist
$\lambda_1,\lambda_2,\lambda_3,\lambda_4\in S$ such that:

$v_{\sigma 2}=\lambda_1 \beta +\lambda_2\gamma$,

$-v_{\sigma 1}=\lambda_3 v'_{\sigma 1}\beta+\lambda_4 v''_{\sigma
1}\gamma $,

$0=\lambda_1 v''_{\sigma 1}+\lambda_3 v_{\sigma 2}$,

$0=\lambda_2 v''_{\sigma 1}+\lambda_4 v_{\sigma 2}$.

It follows that $v_{\sigma 2}|\lambda_1$ and $v_{\sigma
2}|\lambda_2$ and so we get $1\in (\beta,\gamma )$. Contradiction!
If $(v_{\sigma 2}v''_{\sigma 1}, \beta)\not\cong 1$ then we are in
the subcase (II5), (II6), $\ldots$.  In the same way we treat (II2),
(II3), (II4).

\begin{center}{\bf Subcase:}\; $\alpha_3=v'_{\sigma 1} \beta$,
  $\alpha_4=v''_{\sigma 1} \gamma$
\end{center}

We see that the relations given by the columns of the following
matrix

$$\left(%
\begin{array}{ccccc}

v_{\sigma 2}& \beta    &\gamma & 0 & 0    \\
 -v_{\sigma 1} &0 & 0 & \alpha_3 & \alpha_4   \\
  0 &  -v''_{\sigma 1} & 0 &  -v_{\sigma 2} & 0\\
 0 & 0 & -v'_{\sigma 1}  & 0 & -v_{\sigma 2}  \\

\end{array}%
\right),$$

are elements in  $\Omega^1_S(J)$. The columns 2,3 together with
the last two columns divided by $(\beta,v_{\sigma 2})$,
respectively $(\gamma,v_{\sigma 2})$ are part of a minimal system
of generators. Since $\mu\bigl(\Omega^1_S(J)\bigr)=4$ we see that
the first column is a linear combination
 of the others as above. Thus there exist
$\lambda_1,\lambda_2,\lambda_3,\lambda_4\in S$ such that:

$v_{\sigma 2}=\lambda_1 \beta +\lambda_2\gamma$,

$-v_{\sigma 1}=\lambda_3 v'_{\sigma 1}\beta/(\beta,v_{\sigma 2})
+\lambda_4 v''_{\sigma 1}\gamma/(\gamma,v_{\sigma 2}) $,

$0=\lambda_1 v''_{\sigma 1}+\lambda_3 v_{\sigma
2}/(\beta,v_{\sigma 2})$,

$0=\lambda_2 v'_{\sigma 1}+\lambda_4 v_{\sigma
2}/(\gamma,v_{\sigma 2})$.

It follows that $v_{\sigma 2}/(\beta,v_{\sigma 2})| \lambda_1$ and
$v_{\sigma 2}/(\gamma,v_{\sigma 2})|\lambda_2$ and so we get $1\in
(\beta,\gamma )$ which is false as above if $(\beta,v_{\sigma
2})\cong 1$, $(\gamma,v_{\sigma 2})\cong 1$. Clearly
$\beta,\gamma$ cannot be multiples of $v_{\sigma 2}$ because
otherwise $J$ is only 3 generated. Thus we may suppose for example
$\beta=v'_{\sigma 2}$.

Then $J=(v_{\sigma 1},v_{\sigma 2},v'_{\sigma 1}v'_{\sigma 2},
v''_{\sigma 1} \gamma)$ and the matrix factorizations of
$\Omega^2_R\bigl(J/(f)\bigr)$ is given by the following matrices
$A,B$:

$$A=\left(%
\begin{array}{ccccc}
 0 & -\gamma & -v'_{\sigma 2} & w_{\sigma 1} & 0\\
v'_{\sigma 1} & 0 & 0 &  w_{\sigma 2} & -\gamma v''_{\sigma 1}\\
 -v''_{\sigma 2} & 0 &  v''_{\sigma 1} & 0 & 0 \\
0 &  v'_{\sigma 1} & 0 & 0 & v_{\sigma 2}\\
0 & -w_{\sigma 2} & 0 & 0 & w_{\sigma 1}v''_{\sigma 1}\\
\end{array}%
\right),$$

$$B=\left(%
\begin{array}{ccccc}
-w_{\sigma 2} v''_{\sigma 1} & w_{\sigma 1} v''_{\sigma 1} &
-w_{\sigma 2}
v'_{\sigma 2} & 0 &  v''_{\sigma 1}\gamma\\
0 & 0 & 0 & w_{\sigma 1}v''_{\sigma 1} & -v_{\sigma 2}\\
 -w_{\sigma 2}v''_{\sigma 2} & w_{\sigma 1}v''_{\sigma 2} &
w_{\sigma 1}v'_{\sigma 1} & 0 & \gamma v''_{\sigma 2}\\
v_{\sigma 1} & v_{\sigma 2} & v'_{\sigma 1}v'_{\sigma 2} & \gamma
v''_{\sigma 1} & 0\\
0 & 0 & 0& w_{\sigma 2} & v'_{\sigma 1}\\
\end{array}%
\right).$$

We may add to $\gamma$ multiples of $v'_{\sigma 1}$ because this
means to add to $\alpha_4$ multiples of $v_{\sigma 1}$. Also
adding multiples of the  column 4 of $A$ to the column 2 and then
adding multiples of the row 5 to the row 2 we see that the result
is just the  addition of  some multiples of $w_{\sigma 1}$ to
$\gamma$. On the other hand adding some multiples of the row 5 to
the row 1 and then
 adding some multiples of the column 4 to the column 5 we see that the result
is just the  addition of  some multiples of $w_{\sigma 2}$ to
$\gamma$.

 So, after some elementary transformations on $A,$ we may suppose $\gamma$
to be a polynomial in $ v''_{\sigma 2}$
 and it is enough to
see that deg$(\gamma)=1$. However adding to $\alpha_4$ multiples
of $\alpha_2$  we may add to $\gamma$ multiples of $ v_{\sigma
2}$. Since $ v_{\sigma 2}^{\prime\prime 2}\in (v_{\sigma
2},w_{\sigma 2})$ we may suppose  deg$(\gamma)= 1$, that is
$\gamma=qv''_{\sigma 2}$ for a certain nonzero
 constant  $q$.
Now we multiply the row 1 of $A$ with a $q^{-1}$ then the columns
3,4 with $q$, then rows 2,3 with $q^{-1}$ and finally the column 1
with  $q$. So  we reduce to the case $q=1$.

Thus $J=J_{1\sigma}(a,b,u)$. If we take $\beta=v''_{\sigma 2}$
then similarly we get $J=J_{2\sigma}(a,b,u)$. If
$(\gamma,v_{\sigma 2})\not\cong 1$ we get similarly
$J=J_{2\sigma}(a,b,u)$, $J=J_{1\sigma}(a,b,u)$.

{\bf Subcase $\alpha_3=v'_{\sigma 1} \beta$, $\alpha_4=v''_{\sigma
2} \gamma$}

We see that the relations given by the columns of the following
matrix

$$\left(%
\begin{array}{ccccc}

v_{\sigma 2}& \beta    &v''_{\sigma 2}\gamma & 0 & 0    \\
 -v_{\sigma 1} &0 & 0 & \alpha_3 & \gamma   \\
  0 &  -v''_{\sigma 1} & 0 &  -v_{\sigma 2} & 0\\
 0 & 0 & -v_{\sigma 1}  & 0 & -v'_{\sigma 2}  \\

\end{array}%
\right),$$

are elements in  $\Omega^1_S(J)$. As in the above subcase the
columns 2,5 together with the columns 3,4 divided by
$(\gamma,v_{\sigma 2})$, respectively $(\beta,v_{\sigma 2}),$ form
a minimal system of generators in $\Omega^1_S(J)$.
 Thus the first column is a linear combinations
 of the others and there exist
$\lambda_1,\lambda_2,\lambda_3,\lambda_4\in S$ such that:

$v_{\sigma 2}=\lambda_1 \beta +\lambda_2 v''_{\sigma 1}\gamma/
(\gamma,v_{\sigma 1})$,

$-v_{\sigma 1}=\lambda_3 v'_{\sigma 1}\beta/(\beta,v_{\sigma 2})
+\lambda_4 \gamma $,

$0=\lambda_1 v''_{\sigma 1}+\lambda_3 v_{\sigma
2}/(\beta,v_{\sigma 2})$,

$0=\lambda_2 v_{\sigma 1}/(\gamma,v_{\sigma 1})+\lambda_4
v'_{\sigma 2}$.

It follows $v'_{\sigma 2}|\lambda_2 $ and from the first identity
we see that we get $v_{\sigma 2}|\lambda_1\beta$. If $\lambda_1=0$
then $\lambda_3 =0$ and so $\gamma|v_{\sigma 1}$. If $\gamma$ is a
multiple of $v''_{\sigma 1}$ then we are in the preceding subcase.
If $\gamma$ is a multiple of $v'_{\sigma 1}$ we change $\alpha_3$
with  $\alpha_4$ and so we may suppose the new $\beta$ to be a
multiple of $v''_{\sigma 2}$. This is exactly one possibility
which follows from $\lambda_1\not=0$ which we treat now. From
first identity we see that deg$(\lambda_1\beta) =2$ and so either
$\beta=v'_{\sigma 2}$ or $\beta=v''_{\sigma 2}$. The first
situation lead us to the proceeding subcase, that is $J\in {\cal
J}$. So we may  suppose $\beta=v''_{\sigma 2}$.\\
Thus we reduce to the case $J=(v_{\sigma 1},v_{\sigma
2},v'_{\sigma 1}v''_{\sigma 2}, v''_{\sigma 2} \gamma)$. Adding to
$\gamma$ multiples of $v'_{\sigma 1},v'_{\sigma 2}$ this
 means  to add some multiples of $\alpha_3,\alpha_2$ we may suppose
that $\gamma\in (v''_{\sigma 1},v''_{\sigma 2})$. As before
deg$(\alpha_4)=1$ because for instance $(\alpha_2,\alpha_3)$
contains $v''_{\sigma 2} (v'_{\sigma 1},v'_{\sigma 2})^2$. Then
$\gamma=\tau_1 v''_{\sigma 1}+ \tau_2 v''_{\sigma 2}$ for some
$\tau_1,\tau_2\in K$. Then another relation of
 $\Omega^1_S(J)$ is the transpose of
$(\tau_1 v''_{\sigma 2}, 0, \tau_2v''_{\sigma 2},-v'_{\sigma 1})$.
If $\tau_1\not =0$ then this relation together with the last 4
columns of the previous matrix (some of them divided by something)
form 5 elements from a minimal system of generators of
$\Omega^1_S(J)$. Contradiction! Thus $\tau_1 =0$ and $J\in {\cal
T}$.\\ With similar procedures we treat the other cases.
\end{proof}

Set
$$\rho_{\sigma}(a,b,u)=\left(%
\begin{array}{ccccc}
  0 & w_{\sigma 1} & -v'_{\sigma 2} &-x_j & 0\\
v'_{\sigma 1} & w_{\sigma 2} & 0 &0 & -x_jv''_{\sigma 1}\\
- v''_{\sigma 2} & 0 & v''_{\sigma 1} & 0 & 0\\
0 &0 & 0 & v'_{\sigma 1}&  v_{\sigma 2}\\
0 & 0& 0 & -w_{\sigma 2}& w_{\sigma 1}v''_{\sigma 1}\\
\end{array}%
\right),$$

$$\omega_{\sigma}(a,b,u)=\left(%
\begin{array}{ccccc}
 -w_{\sigma 2}v''_{\sigma 1} & w_{\sigma 1}v''_{\sigma 1} &
-w_{\sigma 2}v'_{\sigma 2} & 0 & x_jv''_{\sigma 1}\\
v_{\sigma 1} & v_{\sigma 2}& v'_{\sigma 1}v'_{\sigma 2 } & x_jv''_{\sigma 1} & 0\\
-w_{\sigma 2}v''_{\sigma 2} & -w_{\sigma 1}v''_{\sigma 2} &
w_{\sigma 1}v'_{\sigma 1} & 0 & x_jv_{\sigma 2 }^{\prime\prime}\\
0 & 0 & 0 & w_{\sigma 1}v''_{\sigma 1} & -v_{\sigma 2}\\
0 & 0 & 0 & w_{\sigma 2} & v'_{\sigma 1}\\
\end{array}%
\right),$$

$$\mu_{ \sigma}(a,b,u)=\left(%
\begin{array}{ccccc}
  0 & w_{\sigma 1} & v''_{\sigma 2} & 0 & 0 \\
 - v'_{\sigma 1} &  w_{\sigma 2} & 0 & 0 & x_j \\
 v'_{\sigma 2} & 0 &  - v''_{\sigma 1} &  x_j & 0 \\
  0 & 0 & 0 & - v'_{\sigma 1} & -v'_{\sigma 2} \\
 0 & 0 & 0 & w_{\sigma 2}v''_{\sigma 2} & - v''_{\sigma 1}w_{\sigma 1} \\
\end{array}%
\right),$$

$$\nu_{ \sigma}(a,b,u)=\left(%
\begin{array}{ccccc}
   v''_{\sigma 1}w_{\sigma 2} & - v''_{\sigma 1}w_{\sigma 1}&  v''_{\sigma 2}w_{\sigma 2} & 0 & -x_j \\
 v_{\sigma 1} &  v_{\sigma 2} &  v''_{\sigma 2} v'_{\sigma 1}  & x_j v''_{\sigma 2} & 0\\
v'_{\sigma 2}w_{\sigma 2} & -v'_{\sigma 2}w_{\sigma 1} & -v'_{\sigma 1}w_{\sigma 1} & -x_jw_{\sigma 1}& 0\\
0 & 0 & 0 & -v''_{\sigma 1}w_{\sigma 1}& v'_{\sigma  2}\\
0& 0 &  0 &  -v''_{\sigma 2} w_{\sigma 2} & -v'_{\sigma 1} \\
\end{array}%
\right),$$

and

$$\bar{\mu}_{ \sigma}(a,b,u)=\left(%
\begin{array}{ccccc}
  0 & w_{\sigma 2} & v'_{\sigma 1} & 0 & 0 \\
 - v''_{\sigma 2} &  w_{\sigma 1} & 0 & 0 & x_s \\
 v''_{\sigma 1} & 0 &  - v'_{\sigma 2} &  x_s & 0 \\
  0 & 0 & 0 & - v''_{\sigma 2} & -v''_{\sigma 1} \\
 0 & 0 & 0 & w_{\sigma 1}v'_{\sigma 1} & - v'_{\sigma 2}w_{\sigma 2} \\
\end{array}%
\right),$$

$$\bar{\nu}_{ \sigma}(a,b,u)=\left(%
\begin{array}{ccccc}
v'_{\sigma 2}w_{\sigma 1} & - v'_{\sigma 2}w_{\sigma 2}&  v'_{\sigma 1}w_{\sigma 1} & 0 & -x_s \\
v_{\sigma 2} &  v_{\sigma 1} &  v'_{\sigma 1} v''_{\sigma 2}  & x_sv'_{\sigma 1} & 0\\
v''_{\sigma 1}w_{\sigma 1} & -v''_{\sigma 1}w_{\sigma 2} & -v''_{\sigma 2}w_{\sigma 2} & -x_sw_{\sigma 2}& 0\\
0 & 0 & 0 & -v'_{\sigma 2}w_{\sigma 2}& v''_{\sigma  1}\\
0& 0 &  0 &  -v'_{\sigma 1} w_{\sigma 1} & -v''_{\sigma 2} \\
\end{array}%
\right).$$

\begin{thm}\label{thm5.3}
Let $$\mathcal{E}=\{\Co\bigl(\rho_{\sigma}(a,b,u)\bigr),\
\Co\bigl(\mu_{\sigma}(a,b,u)\bigr),\
\Co\bigl(\bar{\mu}_{\sigma}(a,b,u)\bigr) \ |\ \sigma, a,b,u \}$$
\begin{itemize}
    \item [(i)]  The set $\mathcal{E}$ contains only
 indecomposable, graded, non--orientable, $5$--generated MCM $R$--modules of
 rank $2$.
    \item [(ii)] Every indecomposable, graded,
 non--orientable, $5$--generated MCM module over $R$ of rank $2$
 is isomorphic with one module of $\mathcal{E}.$
    \item [(iii)] All the modules of the set $\mathcal{E}$ are
 non--isomorphic. In particular, there are $162$ isomorphism classes of
 indecomposable, graded, non--orientable MCM modules over $R$ of rank two, with
 $5$ generators.
\end{itemize}
\end{thm}

\begin{proof}
$(i)$.The pairs of  matrices $(\rho_{\sigma}(a,b,u),
\omega_{\sigma}(a,b,u))$ and
$(\mu_{\sigma}(a,b,u),\nu_{\sigma}(a,b,u))$ have been obtained
from the pairs $(\rho_{1 \sigma}(a,b,u), \omega_{1
\sigma}(a,b,u)),$ respectively $(\mu_{1
\sigma}(a,b,u),\nu_{1\sigma}(a,b,u))$ by elementary operations on
rows and columns. The pair  $(\bar{\mu}_{1
\sigma}(a,b,u),\bar{\nu}_{1\sigma}(a,b,u))$ is a matrix
factorization corresponding to $\Omega_R^2(T_{5\sigma}(a,b,u)).$
By the above Lemma, the set $\mathcal{E}$ satisfies the part $(i)$
of the theorem. For the proof of indecomposability we may proceed
as in the last part of the proof of Theorem \ref{thm3.2}. \\
$(ii)$.Preserving the notations of Lemma \ref{lem5.1}, set
$$\mathcal{J}_i=\{J_{i \sigma}(a,b,u)\ |\ \sigma, a,b,u\}, i=1,2$$
and $$\mathcal{T}_i=\{T_{i \sigma}(a,b,u)\ |\ \sigma, a,b,u\},$$
for $1\leq i\leq 8$. We claim that
$$\mathcal{J}_1=\mathcal{J}_2,$$
$$\mathcal{T}_1=\mathcal{T}_2=\mathcal{T}_3=\mathcal{T}_4$$ and
$$\mathcal{T}_5=\mathcal{T}_6=\mathcal{T}_7=\mathcal{T}_8.$$
Indeed, take, for instance, $J_{2 \sigma}(a,b,u)=(v_{\sigma
1},v_{\sigma 2},v'_{\sigma 1}v''_{\sigma 2},v''_{\sigma
1}v'_{\sigma 2})\in \mathcal{J}_2$.Then we may find
$\overline{v}'_{\sigma 1}, \overline{v}'_{\sigma
2},\overline{v}''_{\sigma 1}, \overline{v}''_{\sigma 2},$
depending on some other cubic roots of $-1,$ let us say $n,p,$ and
$v$, a cubic root of unity different from $1$, such that $$J_{2
\sigma}(a,b,u)=J_{1 \sigma}(n,p,v)=(\overline{v}_{\sigma
1},\overline{v}_{\sigma 2}, \overline{v}'_{\sigma 1}
\overline{v}'_{\sigma 2},\overline{v}''_{\sigma 1}
\overline{v}''_{\sigma 2}).$$\\
$(iii)$.In order to check that  the modules of the list are
pairwise non--isomorphic, we have to prove that the matrices of
the set
\[
\mathcal{E}'=\{\rho_{\sigma}(a,b,u),\mu_{\sigma}(a,b,u),\bar{\mu}_{\sigma}(a,b,u)\|\
\sigma,a,b,u\}
\]
are pairwise non--equivalent. As in the proof of Theorem
\ref{prop4.6}, if $A,B\in \mathcal{E}'$ are two equivalent matrices, then the
matrices $\widetilde{A}$ and $\widetilde{B},$ obtained by reducing the entries
of $A,$ respectively $B,$ modulo $\mathfrak{m}^2,$ are also equivalent by some
scalar and invertible matrices. We observe that the "reduced" matrix
$\tilde{\rho}_{\sigma}(a,b,u)$ has the entries of the last column zero and the
matrices $\tilde{\mu}_{\sigma}(a,b,u), \tilde{\bar{\mu}}_{\sigma}(a,b,u)$ have
the entries of the last row zero. If $\tilde{\rho}_{\sigma}(a,b,u)\sim
\tilde{\mu}_{\tau}(n,p,v),$ for some $\sigma, a,b,u,\tau,n,p,v,$ then there
exist some invertible scalar $5\times 5$ matrices $U,V$ such that $$U\cdot
\tilde{\rho}_{\sigma}(a,b,u)=\tilde{\mu}_{\tau}(n,p,v)\cdot V.$$
Looking at
the last column in this equality, we obtain that $V$ must have the last column
zero, contradiction. In the same way we obtain that
$\tilde{\rho}_{\sigma}(a,b,u)\not\sim \tilde{\bar{\mu}}_{\tau}(n,p,v)$.\\ Let
us suppose now that $\tilde{\mu}_{\sigma}(a,b,u)\sim
\tilde{\bar{\mu}}_{\tau}(n,p,v),$ for some $\sigma, a,b,u,\tau,n,p,v,$ and let
$U,V\in \mathcal{M}_{5 \times 5}(K)$ invertible such that
\[
\tilde{\mu}_{\sigma}(a,b,u)\cdot U=V\cdot \tilde{\bar{\mu}}_{\tau}(n,p,v).
\]
We compare the entries of the fourth column in the above equality. Let
$\tau=(e\ f\ t)$.  For $t\in \{i,j\},$ which implies $\sigma\neq \tau,$ we
obtain that all the entries of the fourth column in $U$ are zero,
contradicting $U$ invertible. If $t\not\in \{i,j\},$ which implies
$\sigma=\tau$ and $t=s,$ we obtain that all the entries of the third column in
$V$ are zero, contradiction.\\ In the same way we may prove that if
$\tilde{\mu}_{\sigma}(a,b,u)\sim \tilde{\mu}_{\tau}(n,p,v)$ or
$\tilde{\bar{\mu}}_{\sigma}(a,b,u)\sim \tilde{\bar{\mu}}_{\tau}(n,p,v)$ or
$\tilde{\rho}_{\sigma}(a,b,u)\sim \tilde{\rho}_{\tau}(n,p,v),$ then
$(\sigma,a,b,u)=(\tau,n,p,v).$
\end{proof}

\begin{cor}\label{cor5.4}
Let $$\mathcal{F}=\{\Co\bigl(\omega_{\sigma}(a,b,u)\bigr),\
\Co\bigl(\nu_{\sigma}(a,b,u)\bigr),\
\Co\bigl(\bar{\nu}_{\sigma}(a,b,u)\bigr) \ |\ \sigma, a,b,u \}$$
\begin{itemize}
    \item [(i)]  The set $\mathcal{F}$ contains only
 indecomposable, graded, non--orientable, $5$--generated MCM $R$--modules of rank $3$.
    \item [(ii)] Every indecomposable, graded,
 non--orientable, $5$--generated MCM module over $R$ of rank $3$
 is isomorphic with one module of $\mathcal{F}.$
    \item [(iii)] All the modules of the set $\mathcal{F}$ are
 non--isomorphic. In particular, there are $162$ isomorphism classes of
 indecomposable, graded, non--orientable MCM modules over $R$ of rank three, with
 $5$ generators.
\end{itemize}
\end{cor}

\begin{proof}
The map $M\mapsto \Omega_R^1(M)$ is a bijection between the
$5$--generated, indecomposable, graded, MCM $R$--modules of rank
$2$ and the $5$--generated, indecomposable, graded, MCM
$R$--modules of rank $3.$
\end{proof}

\begin{lem}\label{lem5.5} There exist no graded, indecomposable,
non--orientable, rank two, six--generated MCM modules.
\end{lem}

\begin{proof} Suppose there exist such MCM module $M$. Then $M\cong
  \Omega^2_R \bigl(J/(f)\bigr)$ for a certain 5--generated ideal
  $J=(\alpha_1,\alpha_2,\alpha_3, \alpha_4,\alpha_5)$ of $S$ as hinted at in
  the first part of Section 4. Then any 4 elements from the $\alpha_t$ must
  generate an ideal $J''$ in ${\cal J}\cup {\cal T}$ because otherwise
  $\mu(\Omega^1_S\bigl(J''/(f)\bigr)>4$ and so obvious
  $\mu(\Omega^1_S(J/(f))>5$. So we may suppose $\alpha_t=v_{\sigma t}$ for
  $t=1,2$ and after some permutations $\alpha_3=v'_{\sigma 1}v''_{\sigma 2}$.
  Set $J'=(\alpha_1,\alpha_2,\alpha_3)$. If $(J',\alpha_4) \in {\cal J}$ then
  $(J',\alpha_5)\not\in {\cal J}$ because otherwise we get
  $\alpha_4=\alpha_5$. Thus $(J',\alpha_5)\in {\cal T}$ and so either
  $\alpha_5= v_{\sigma 2}^{\prime\prime 2}$ or $\alpha_5= v_{\sigma 1}^{\prime
    2}$.  But then $(\alpha_1,\alpha_2,\alpha_4,\alpha_5)\not \in {\cal J}\cup
  {\cal T}$.  If $(J',\alpha_4)\not\in {\cal J}$ and $(J',\alpha_5)\not\in
  {\cal J}$ then $(J',\alpha_4),(J',\alpha_5)\in {\cal T}$ and so $\alpha_4=
  v_{\sigma 2}^{\prime\prime 2}$ and $\alpha_5= v_{\sigma 1}^{\prime 2}$ or
  conversely. But then $(\alpha_1,\alpha_2,\alpha_4,\alpha_5)\not \in {\cal
    J}\cup {\cal T}$.
\end{proof}

\begin{cor}\label{cor5.6}
  There exist no indecomposable, graded, non--orientable, rank 4, 6--generated
  MCM modules.
\end{cor}


\section{Orientable, rank 2, 6--generated MCM modules}

Let $S = K[x_1, x_2, x_3, x_4]$, and $R = S/(f)$, $f = x_1^3 + x_2^3 + x_3^3 +
x_4^3$.

We have proved that a non--free graded orientable six--generated
MCM $R$--module corresponds to a skew symmetric homogeneous matrix
over $S$ of order 6, whose determinant is $f^2$.

Let $\Lambda$ be such a matrix.   Notice that $\Lambda$ has linear
entries and the matrix $\underline{\Lambda} := \Lambda|_{x_4 =
0}$, obtained from $\Lambda$ by restricting the entries to $x_4 =
0$, is a homogeneous matrix over $S_3 = K[x_1, x_2, x_3]$, whose
determinant is $f_3^2$, where $f_3 = x_1^3 + x_2^3 + x_3^3$.
Therefore, Coker$\underline{\Lambda}$ defines a graded rank two,
six--generated MCM over $R_3 = S_3/(f_3)$.   These modules were
explicitly described in \cite{LPP}.

\begin{lem}\label{lemma6.1}
  Let $M$ be a non--free graded orientable six--generated MCM module over
  $R$.   Then the restriction of $M$ to the curve defined by $f = x_4 = 0$
  splits into a direct sum of a 3---generated MCM of rank 1 and its dual.
  Especially, there exists $\lambda \in V(f_3)\smallsetminus\{P_0\}$ and a
  skew symmetric matrix $\Gamma \in \km_{6\times6}(K)$, such that $M$ is the
  cokernel of a map given by the matrix $\Lambda = x_4 \cdot \Gamma + \left(
    \begin{smallmatrix}
      0 & -\alpha^t_\lambda\\ \alpha_\lambda & \phantom{-}0
    \end{smallmatrix}\right)$.

(The same notations as in \cite{LPP} and in ``Preliminaries''.)
\end{lem}

\begin{proof}
  Let $\Lambda_1$ be a skew symmetric homogeneous matrix over $S$,
  corresponding to $M$, and denote $\underline{\Lambda_1} =
  \Lambda_1|_{x_4=0}$.  Suppose that the MCM $S_3$--module corresponding to
  $\underline{\Lambda_1}$ is indecomposable.  Then we can generate it as
  described in Theorem 4.2 and Lemma 5.4 from \cite{LPP}.  Denote with $D$ the
  matrix which we obtain by this means.

  Since $D \sim \underline{\Lambda_1}$, and $\underline{\Lambda_1}$ is skew
  symmetric, there exist two invertible matrices $U, V \in \km_{6\times 6}(K)$
  such that $U \cdot D \cdot V + (U \cdot D \cdot V)^t = 0$.  Therefore, there
  exists $T \in \km_{6 \times 6}(K)$ an invertible matrix such that $T \cdot D
  + (TD)^t = 0$.  (Take $T = (V^t)^{-1} \cdot U$.)

With the help of \textsc{Singular}, we find that, in fact, there
is no invertible matrix $T$  such that $T \cdot D$ is skew
symmetric.   Therefore, the module corresponding to
$\underline{\Lambda_1}$ should decompose.

\begin{verbatim}
//First, we generate the matrix D

LIB"matrix.lib"; option(redSB); proc reflexivHull(matrix M) {
  module N=mres(transpose(M),3)[3];
  N=prune(transpose(N));
  return(matrix(N));
}

proc tensorCM(matrix Phi, matrix Psi) {
   int s=nrows(Phi);
   int q=nrows(Psi);
   matrix A=tensor(unitmat(s),Psi);
   matrix B=tensor(Phi,unitmat(q));
   matrix R=concat(A,B,U);
   return(reflexivHull(R));
}

proc M2(ideal I) {
   matrix A=syz(transpose(mres(I,3)[3]));
   return(transpose(A));
}

ring R=0,(x(1..3)),(c,dp); qring S=std(x(1)^3+x(2)^3+x(3)^3);
ideal I=maxideal(1); matrix C=M2(I);

ring R1=(0,a),(x(1..3),e,b),lp; ideal I=x(1)^3+x(2)^3+x(3)^3,
         (a-1)^3+b3+1,e*b+a2-3*a+3,e*a-b2;
qring S1=std(I);

matrix B[3][3]=       0,  x(1)-(a-1)*x(3),    x(2)-b*x(3),
               x(1)+x(3),  -x(2)-x(3)*b,    -x(3)*e,
                    x(2),        x(3)*e,  -x(1)+(-a+2)*x(3);
matrix C=imap(S,C); matrix D=tensorCM(C,B);

//We check the existence of the invertible matrix T

ring R2=0,(x(1..3),a,e,b,t(1..36)),dp; ideal
I=x(1)^3+x(2)^3+x(3)^3,(a-1)^3+b3+1,e*b+a2-3*a+3,e*a-b2; qring
S2=std(I); matrix D=imap(S1,D); matrix T[6][6]=t(1..36); matrix
A=T*D+transpose(T*D); ideal I=flatten(A); ideal
I1=transpose(coeffs(I,x(1)))[2]; ideal
I2=transpose(coeffs(I,x(2)))[2]; ideal
I3=transpose(coeffs(I,x(3)))[2]; ideal J=I1+I2+I3+ideal(det(T)-1);
ideal L=std(J);
 L;

L[1]=1
\end{verbatim}
\texttt{//Therefore, there does not exist an invertible matrix $T$ such that
  $T\cdot D$ skew symmetric.}

So, after some linear transformations, $\underline{\Lambda_1}$
decomposes into two matrices of order three and rank one with
determinant $f_3 = x_1^3 + x_2^3 + x_3^3$, which correspond to two
points $\lambda_1, \lambda_2$ in $V(f_3)\smallsetminus\{P_0\}$,
$P_0 = [-1:0:1]$.   Let us denote them by $A$ and $B$.   We can
consider $A = \alpha_{\lambda_1}$, $B = \alpha_{\lambda_2}$.

Since $\underline{\Lambda_1}$ is skew symmetric, there exists an
invertible matrix $U \in \km_{6 \times 6}(K)$ such that $U \cdot
\left(
  \begin{smallmatrix}
    A & 0\\0 & B
  \end{smallmatrix}\right)$ is skew symmetric.   Therefore, if we consider $U
= \left(
  \begin{smallmatrix}
    U_1 & U_2\\U_3 & U_4
  \end{smallmatrix}\right)$, we have the following equalities:
\[
\left\{
  \begin{array}{ccc}
U_1 \cdot A + (U_1 \cdot A)^t & = & 0\\
U_4 \cdot B + (U_4 \cdot B)^t & = & 0\\
U_2 \cdot B + \phantom{(}A^t \cdot U_3^t\phantom{)} & = & 0\\
U_3 \cdot A + \phantom{(}B^t \cdot U^t_2\phantom{)} & = & 0
\end{array}\right.
\]
So $U_1 \cdot \alpha_{\lambda_1}$ and $U_4 \cdot
\alpha_{\lambda_2}$ are skew symmetric, so they have only zeros on
the main diagonal.   Since the entries of the second and third
line and column of $\alpha_{\lambda_1}$ and $\alpha_{\lambda_2}$
are linearly independent, we easily obtain that $U_1 = U_4 = 0$.
Therefore, $U_2$ and $U_3$ are invertible matrices and $B = -
U_2^{-1} \cdot A^t \cdot U_3^t$.

We have obtained $\underline{\Lambda_1} \sim \left(
  \begin{smallmatrix}
    \alpha_{\lambda_1} & 0\\0 & \alpha_{\lambda_1}^t
  \end{smallmatrix}\right) \sim \left(
  \begin{smallmatrix}
    0 & -\alpha^t_{\lambda_1}\\ \alpha_{\lambda_1} & \phantom{-}0
  \end{smallmatrix}\right)$.

Therefore, there exists $\Gamma \in \km_{6 \times 6}(K)$ skew
symmetric and $\lambda \in V(f_3)\smallsetminus\{P_0\}$ such that
$\Lambda_1 \sim \Lambda = x_4 \cdot \Gamma + \left(
  \begin{smallmatrix}
    0 & -\alpha^t_{\lambda}\\ \alpha_{\lambda} & \phantom{-}0
  \end{smallmatrix}\right)$.   We can write $\Gamma = \left(
  \begin{smallmatrix}
    \Gamma_1 & -\Gamma_2^t\\\Gamma_2 & \phantom{-}\Gamma_3
  \end{smallmatrix}\right)$, $\Gamma_i \in \km_{3\times 3}(K)$, $i =
1,2,3$, $\Gamma_1$ and $\Gamma_3$ skew symmetric.
\end{proof}

\begin{rem}[Notation]\label{remark6.2}
  For any $\lambda = [a : b :c] \in V(f_3)\smallsetminus\{P_0\}$ there exists
  a unique point in $V(f_3)\smallsetminus\{P_0\}$ which we denote as
  $\lambda^t$, such that $\alpha_\lambda^t \sim \alpha_{\lambda^t}$.   We find
  $\lambda^t = [c:b:a]$.\footnote{If $\lambda$ corresponds to the 3--generated
  rank 1 MCM $N$, then $\lambda^t$ corresponds to its dual $N^\vee$.}

For $\lambda = [a:b:1]$ we denote with $U_\lambda$ and $V_\lambda$
two invertible matrices such that $U_\lambda \cdot
\alpha_\lambda^t = \alpha_{\lambda^t} \cdot V_\lambda$.

If $a \not= 0$, then we can take $U_\lambda = \left(
  \begin{smallmatrix}
    \phantom{-}b^2 && b(a+1) && -(a+1)^2\\
    -(a+1)^2 && b^2 && -b(a+1)\\
\phantom{-}b(a+1) && (a+1)^2 && b^2
  \end{smallmatrix}\right)$ and $V_\lambda = U_\lambda^t$.

If $a = 0$, then we can take $U_\lambda = \left(
  \begin{smallmatrix}
    -b^2 && -b & &1\\
-2b && 1 && b^2\\
\phantom{-} 2b^2 &&2b&&1
  \end{smallmatrix}\right)$ and $V_\lambda = \left(
  \begin{smallmatrix}
    \phantom{-}1 && -2b &&\phantom{-}2b^2\\
-b && -b^2 && -1\\
-b && -b^2 &&\phantom{-}2
  \end{smallmatrix}\right)$.

Notice that for $\lambda = [1:b:1] \in V(f_3)$, $\lambda^t =
\lambda$ and for all other $\lambda \in
V(f_3)\smallsetminus\{P_0\}$, $\lambda \not= \lambda^t$.\qed
\end{rem}

\begin{rem}\label{remark6.3}
  For any $\lambda = [1:b:0] \in V(f_3)\smallsetminus\{P_0\}$ and any $\Lambda
  = x_4 \cdot \Gamma + \left(
    \begin{smallmatrix}
      0 & -\alpha^t_\lambda\\
\alpha_\lambda & \phantom{-}0
    \end{smallmatrix}\right)$ skew symmetric with $\det \Lambda = f^2$, we
  have  $\lambda^t = [0 :b : 1]$ in $V(f_3)\smallsetminus \{P_0\}$ and
  $\Lambda^\prime = x_4 \Gamma^\prime + \left(
    \begin{smallmatrix}
      0 & - \alpha_{\lambda^t}^t\\\alpha_{\lambda^t} & 0
    \end{smallmatrix}\right)$ skew symmetric with
  $\det \Lambda^\prime = f^2$ such that $\Lambda \sim \Lambda^\prime$.

Indeed, take $\Lambda^\prime = U \cdot \Lambda \cdot U^t$ where $U
= \left(
  \begin{smallmatrix}
    0 & T_1\\T_2 & 0
  \end{smallmatrix}\right)$, $T_1 = \frac{1}{\sqrt{3}} \left(
  \begin{smallmatrix}
    b & b^2 & 0\\
\phantom{-}0 & \phantom{-}1 & -b^2\\
\phantom{-}2 & \phantom{-}0 & 1
  \end{smallmatrix}\right)$ and $T_2 = \frac{1}{\sqrt{3}} \left(
  \begin{smallmatrix}
           -1 & b & \phantom{-}b\\
 \phantom{-}2 & b^2 & \phantom{-}b^2\\
         -2b^2 & 1 & -2
  \end{smallmatrix}\right)$.

Therefore, Coker $\Lambda$ and Coker $\Lambda^\prime$ define two
isomorphic MCM modules.   This is the reason why, from now on, we
may only consider the case $\lambda = [a:b:1] \in
V(f_3)\smallsetminus\{P_0\}$.\qed
\end{rem}

\begin{rem}\label{remark6.4}
  Consider $\lambda = [a:b:1] \in V(f_3)\smallsetminus\{P_0\}$ and $\Lambda =
  x_4 \cdot \Gamma + \left(
    \begin{smallmatrix}
      0 & -\alpha_\lambda^t\\
\alpha_\lambda & \phantom{-}0
    \end{smallmatrix}\right)$ as in Lemma \ref{lemma6.1}.   Then there exists
  $\overline{\Lambda} = x_4 \cdot \overline{\Gamma} + \left(
    \begin{smallmatrix}
      \alpha_\lambda & 0\\0 & \alpha_{\lambda^t}
    \end{smallmatrix}\right)$ with $\det \overline{\Lambda} = f^2$ such that
  $\overline{\Lambda} \sim \Lambda$.

Indeed, consider $\overline{\Lambda} = \left(
  \begin{smallmatrix}
    \phantom{-}0 & \id\\-U_\lambda & 0
  \end{smallmatrix}\right) \cdot \Lambda \cdot \left(
  \begin{smallmatrix}
    \id & 0\\
0 & V_\lambda^{-1}
  \end{smallmatrix}\right)$.

We obtain $\overline{\Gamma} = \left(
  \begin{smallmatrix}
    \phantom{-}\Gamma_2 && \Gamma_3 \cdot V_\lambda^{-1}\\
-U_\lambda \cdot \Gamma_1 && U_\lambda \cdot \Gamma_2^t \cdot
V_\lambda^{-1}
  \end{smallmatrix}\right)$.\qed
\end{rem}

\begin{lem}\label{lemma6.5}
  Consider $\Lambda = x_4 \cdot \Gamma + \left(
    \begin{smallmatrix}
      0 & -\alpha_\lambda^t\\
\alpha_\lambda & \phantom{-}0
    \end{smallmatrix}\right)$ as above.   Then the MCM module $M$
  corresponding to $\Lambda$ is indecomposable if and only if $\Gamma_1 \not=
  0$ or $\Gamma_3 \not= 0$.
\end{lem}

\begin{proof}
  Suppose $M$ is indecomposable.   If $\Gamma_1 = \Gamma_3 = 0$, then $\left(
    \begin{smallmatrix}
      0 & \id\\\id & 0
    \end{smallmatrix}\right) \cdot \Lambda = \left(
    \begin{smallmatrix}
      x_4 \cdot \Gamma_2 + \alpha_\lambda & & \phantom{-}0\\
0 & & -x_4 \cdot \Gamma_2^t - \alpha_\lambda^t
    \end{smallmatrix}\right)$, so $\Lambda$ decomposes after some linear
  transformation.

This contradicts the indecomposability of $M = \text{ Coker }
\Lambda$, so we must have $\Gamma_1 \not= 0$ or $\Gamma_3 \not=
0$.

Now, let us suppose $\Gamma_1 \not= 0$ or $\Gamma_3 \not= 0$ and
prove that $M$ is indecomposable.

Suppose $M$ decomposes.   Then there exists a matrix $\left(
  \begin{smallmatrix}
    T_1 & 0\\ 0 & T_2
  \end{smallmatrix}\right)$ equivalent to $\Lambda$ with $T_1, T_2$ two
matrices of order three and rank one, with $\det T_1 = \det T_2 =
f$ and $T_1|_{x_4=0} = \alpha_{\lambda_1}$, $T_2|_{x_4=0} =
\alpha_{\lambda_2}$, where $\lambda_1, \lambda_2 \in
V(f_3)\smallsetminus \{P_0\}$.

Since $\Lambda$ is skew symmetric, after some linear
transformations, $\left(
  \begin{smallmatrix}
    \alpha_{\lambda_1} & 0\\ 0 & \alpha_{\lambda_2}
  \end{smallmatrix}\right)$ should also become skew symmetric.   As we saw in
the proof of Lemma \ref{lemma6.1}, this gives $\alpha_{\lambda_2}
\sim \alpha_{\lambda_1}^t$, so $\lambda_2 = \lambda_1^t$.

Using Remark \ref{remark6.3}, there exist $U, V \in \km_{6 \times
6} (K)$ invertible matrices such that $U \cdot \overline{\Lambda}
\cdot V = \left(
  \begin{smallmatrix}
    T_1 & 0\\ 0 & T_2
  \end{smallmatrix}\right) = x_4 \cdot \left(
  \begin{smallmatrix}
    N_1 & 0\\ 0 & N_2
  \end{smallmatrix}\right) + \left(
  \begin{smallmatrix}
    \alpha_{\lambda_1} & 0\\ 0 & \alpha_{\lambda_1^t}
  \end{smallmatrix}\right)$.

Therefore,
\[
\left\{\begin{array}{lcll} U \cdot \left(\begin{array}{cc}
\alpha_\lambda & 0\\ 0 &
    \alpha_{\lambda^t}\end{array}\right) & =
&\left(\begin{array}{cc}\alpha_{\lambda_1} & 0\\ 0 &
    \alpha_{\lambda_1^t}\end{array}\right) & \cdot V^{-1} \qquad (1)\\[1.0ex]
U \cdot \left(\begin{array}{ccc} \Gamma_2 & & \Gamma_3 \cdot
    V_\lambda^{-1}\\-U_\lambda \cdot \Gamma_1 & & U_\lambda \cdot \Gamma_2^t
    \cdot V_\lambda^{-1}\end{array}\right) & = & \left(
  \begin{array}{cc}
N_1 & 0\\0 & N_2
  \end{array}\right) & \cdot V^{-1} \qquad (2)\,.
\end{array}
\right.
\]

Let us consider $U = \left(
  \begin{smallmatrix}
    U_1 & U_2\\U_3 & U_4
  \end{smallmatrix}\right)$ and $V^{-1} = \left(
  \begin{smallmatrix}
    V_1 & V_2\\ V_3 & V_4
  \end{smallmatrix}\right)$ with $U_i, V_i \in \km_{3 \times 3} (K)$, $i =
\overline{1,4}$.

The first system of equations gives:
\[
\left\{\begin{array}{lcl}
U_1 \cdot \alpha_\lambda & = & \alpha_{\lambda_1} \,\cdot V_1\\
U_2 \cdot \alpha_{\lambda^t} & = & \alpha_{\lambda_1} \, \cdot V_2\\
U_3 \cdot \alpha_\lambda & = & \alpha_{\lambda_1^t} \cdot V_3\\
U_4 \cdot \alpha_{\lambda^t} & = & \alpha_{\lambda_1^t} \cdot V_4
\end{array}
\right.
\]
By comparing the coefficients of $x_1, x_2, x_3$ on the left--hand
side and right--hand side of the above equalities, we obtain
easily:
\[
U_i = V_i = K_i \cdot \id_3 \text{ with } K_i \in K,\; i =
\overline{1,4}\,.
\]
Moreover, if $\lambda \not= \lambda_1$, then $K_1 = K_4 = 0$ and
if $\lambda \not= \lambda_1^t$, then $K_2 = K_3 = 0$.   Since $U$
is invertible, we have $\lambda = \lambda_1$ or $\lambda =
\lambda_1^t$.

We know that $\alpha_{\lambda_1} = T_1|_{x_4=0}$ where $T_1$ is a
matrix of order three over $S = K[x_1, x_2, x_3, x_4]$ of rank one
and with determinant $f$.   So Coker $T_1$ is a graded
three--generated rank one MCM $R$--module.   In \cite{EP}, all the
isomorphism classes of such modules are given explicitly.   We
obtain $\alpha_{\lambda_1} \sim \alpha|_{x_4=0}$
or $\alpha_{\lambda_1} \sim \alpha^t|_{x_4=0}$ or $\alpha_{\lambda_1}
\sim \eta|_{x_4=0}$ or $\alpha_{\lambda_1} \sim \nu|_{x_4=0}$.

With the help of computers, we obtain that none of the above
matrices is equivalent to $\alpha_{[1:\ell:1]}$, therefore,
$\lambda_1 \not= \lambda_1^t$.

\begin{alltt}
LIB"matrix.lib"; option(redSB);

ring r=0,(x(1..3),l,a,b,c,d,e,v(1..9),u(1..9)),dp; ideal
I=x(1)^3+x(2)^3+x(3)^3,
        l^3+2,
        a3+1,b3+1,c3+1,d3+1,e2+e+1,bcd-e*a;
qring s=std(I);

proc isomorf(matrix X,matrix Y) {
   matrix U[3][3]=u(1..9);
   matrix V[3][3]=v(1..9);
   matrix C=U*X-Y*V;
   ideal I=flatten(C);
   ideal I1=transpose(coeffs(I,x(1)))[2];
   ideal I2=transpose(coeffs(I,x(2)))[2];
   ideal I3=transpose(coeffs(I,x(3)))[2];
   ideal J=I1+I2+I3+ideal(det(U)-1,det(V)-1);
   ideal L=std(J);
   L;
}


matrix A[3][3]=0,    x(1)-x(3),  x(2)-l*x(3),
       x(1)+x(3), -x(2)-l*x(3), -1/2*l^2*x(3),
            x(2),  1/2*l^2*x(3),        -x(1);

//This is the matrix corresponding to the point (1:l:1)

//We now write the matrices corresponding to the rank one three-
generated MCM modules, restricted to x(4)=0


matrix alpha[3][3]=0,                x(1), -x(3)*b+x(2),
         -x(2)*c+x(1),           -x(3)*b^2, x(3)*b^2*c^2,
                 x(3), x(3)*b*c^2+x(2)*c^2, -x(2)*c-x(1);


matrix alphat=transpose(alpha);

matrix eta[3][3]=0,x(1)+x(2),          x(3),
       x(1)+e*x(2),    -x(3),             0,
              x(3),        0,-x(1)-e^2*x(2);


matrix nu[3][3]=0,x(1)+x(3),           x(2),
    x(1)-a^2*b*x(3),    -x(2),              0,
               x(2),        0,-x(1)+a*b^2*x(3);

isomorf(alpha,A); L[1]=1 isomorf(alphat,A); L[1]=1 isomorf(eta,A);
L[1]=1 isomorf(teta,A); L[1]=1

// Therefore none is isomorphic to \(\alpha\sb{[1:\ell:1]}\) and this means \(\lambda\sb{1} \not= \lambda\sp{t}\sb{1}\).
\end{alltt}

If $\lambda = \lambda_1 \not= \lambda_1^t$ as a solution of the
system (1), we obtain: $U = V = \left(
  \begin{smallmatrix}
K_1 \cdot \id & 0\\ 0 & K_4 \cdot \id
  \end{smallmatrix}\right),\; K_1 \cdot K_4 \not= 0$.

Replacing $U$ and $V$ in (2), we obtain: $\begin{cases}
 K_1 \cdot \Gamma_3 \cdot V_\lambda = 0\\
K_4 \cdot U_\lambda \cdot \Gamma_1 = 0.
\end{cases}$

Since $K_1 \not= 0$, $K_4 \not= 0$ and $U_\lambda, V_\lambda$ are
invertible matrices, we obtain $\Gamma_1 = \Gamma_3 = 0$, which is
a contradiction to our hypothesis.

If $\lambda = \lambda_1^t \not= \lambda$, we obtain as a solution
of (1): $U = V = \left(
  \begin{smallmatrix}
0 & K_2 \cdot \id\\  K_3  \id & 0
  \end{smallmatrix}\right),\; K_2 \cdot K_3 \not= 0$.

Replacing $U$ and $V$ in (2), we obtain: $\begin{cases}
 K_2 \cdot U_\lambda \cdot \Gamma_1 = 0\\
K_3 \cdot \Gamma_3 \cdot V_\lambda = 0.
\end{cases}$
Therefore, we must have again $\Gamma_1 = \Gamma_3 = 0$.
\end{proof}

For each $\lambda = [a:b:1] \in V(f_3)\smallsetminus\{P_0\}$, we
define a family of skew symmetric homogeneous indecomposable
matrices of order six over $S = K[x_1, x_2, x_3, x_4]$ with
determinant $f^2$:
\begin{align*}
\km_\lambda := \left\{\Lambda_{(\lambda,\Gamma)} = x_4 \cdot
\Gamma +\right.
  & \left(\begin{array}{cc}0 & -\alpha_\lambda^t\\\alpha_\lambda &
  0\end{array}\right),\; \det \Lambda_{(\lambda,\Gamma)}= f^2,\\
  \Gamma =
  &\left.\left(\begin{array}{cc} \Gamma_1 & - \Gamma_2^t\\\Gamma_2 &
  \Gamma_3\end{array}\right),\; \Gamma_1, \Gamma_3 \text{ skew symmetric,}\;
  \Gamma_1 \not= 0 \text{ or } \Gamma_3 \not= 0\right\}\,.
\end{align*}
Notice that, as in the proof of Lemma \ref{lemma6.5}, if
$\Lambda_{(\lambda,\Gamma)} \sim \Lambda_{(\lambda^\prime,
\Gamma^\prime)}$, then $\lambda^\prime = \lambda$ or
$\lambda^\prime = \lambda^t$.

\begin{lem}\label{lem6.6}
  Let $\lambda = [a:b:1] \in V(f_3) \smallsetminus \{P_0\}$ with $a \not= 1$.

  \begin{enumerate}
  \item Inside the family $\km_{\lambda}$, two matrices, $\Lambda$ and
  $\Lambda^\prime$, are equivalent if and only if there exists $k \in K^\ast$
  such that $\Lambda^\prime = U_k \cdot \Lambda \cdot U_k^t$, $U_k = \left(
    \begin{smallmatrix}
      k\id & 0\\0 & \frac{1}{k}\id
    \end{smallmatrix}\right)$.   This condition means:
  $\begin{cases}
    \Gamma_2^\prime = \Gamma_2\\
\Gamma_1^\prime = k^2 \cdot \Gamma_1\\
\Gamma_3^\prime = \frac{1}{k^2} \cdot \Gamma_3.
  \end{cases}$

\item A matrix $\Lambda$ from $\km_\lambda$ is equivalent to a matrix
  $\Lambda^\prime$ from $\km_{\lambda^\prime}$, $\lambda^\prime \not= \lambda$
  if and only if $\lambda^\prime = [1:b:a]$
  and $\Lambda^\prime = U_k \cdot \Lambda \cdot U_k^t$, where $k \in K^\ast$
  and $U_k = \left(
    \begin{smallmatrix}
      0 & k \cdot U_\lambda^{-1}\\ -\frac{1}{k} U_\lambda & 0
    \end{smallmatrix}\right)$.
  \end{enumerate}
\end{lem}

\begin{proof}
  We assume $a \not= 0$.  The case $a = 0$ is treated similarly.  Two
  matrices, $\Lambda = \Lambda_{(\lambda, \Gamma)}$ and $\Lambda^\prime =
  \Lambda_{(\lambda^\prime,\Gamma^\prime)}$, are equivalent if and only if
  $\overline{\Lambda}$ and $\overline{\Lambda}^\prime$ are equivalent (see
  Remark \ref{remark6.4}).

If $U$ and $V$ are two invertible matrices such that $U \cdot
\overline{\Lambda} = \overline{\Lambda}^\prime \cdot V$, as in the
proof of Lemma \ref{lemma6.5}, we obtain
\[
\begin{array}{ll}
U = V = \left(\begin{array}{cc}K_1\id & K_2\id\\ K_3\id &
    K_4\id\end{array}\right) \text{ with}& K_1 = K_4 = 0 \text{ if } \lambda
    \not= \lambda^\prime \text{ and}\\
&  K_2 = K_3 = 0 \text{ if } \lambda^\prime
    \not= \lambda^t\,.\end{array}
\]
Since $U \cdot \overline{\Lambda} = \overline{\Lambda}^\prime
\cdot V$, we have:
\begin{equation}
\left(
  \begin{array}{cc}
0 & \id\\ -U_{\lambda^\prime} & 0
  \end{array}\right)^{-1} \cdot U \cdot \left(
  \begin{array}{cc}
0 & \id\\ -U_{\lambda} & 0
  \end{array}\right) \cdot \Lambda \cdot \left(
  \begin{array}{cc}
\id & 0\\ 0 & V_{\lambda}^{-1}
  \end{array}\right) \cdot U^{-1} \cdot  \left(
  \begin{array}{cc}
\id & 0\\ 0 & V_{\lambda^\prime}^{-1}
  \end{array}\right)^{-1} = \Lambda^\prime.  \tag{$\ast$}
\end{equation}

\begin{enumerate}
\item If $\lambda = \lambda^\prime$ then $\lambda^\prime \not=
\lambda^t$, so
  $U = \left(
    \begin{smallmatrix}
      K_1\id & 0\\  0 & K_4\id
    \end{smallmatrix}\right)$ with $K_1 \not= 0$, $K_4 \not= 0$.   So
  ($\ast$) implies: $\left(
    \begin{smallmatrix}
      K_4 \cdot \id & 0\\ 0 & K_1\id
    \end{smallmatrix}\right) \cdot \Lambda \cdot \left(
    \begin{smallmatrix}
      \frac{1}{K_1} \id & 0\\ 0 & \frac{1}{K_4} \id
    \end{smallmatrix}\right) = \Lambda^\prime$.   For $k =
  \sqrt{\frac{K_4}{K_1}}$ and $U_k = \left(\begin{smallmatrix}k\id & 0\\0 &
      \frac{1}{k} \cdot \id\end{smallmatrix}\right)$ we have $\Lambda^\prime =
  U_k \cdot \Lambda \cdot U_k^t$.

\item If $\lambda^\prime = \lambda^t$ then $\lambda^\prime \not=
\lambda$, so
  $U = \left(
    \begin{smallmatrix}
      0 & K_2\id\\ K_3\id & 0
    \end{smallmatrix}\right)$, $K_2 \not= 0$, $K_3 \not= 0$.   Replacing $U$
  in ($\ast$) we obtain:
\[
\Lambda^\prime = \left(
  \begin{array}{cc}
0 & -K_3 U_{\lambda^t}^{-1}\\ -K_2 U_\lambda & 0
  \end{array}\right) \cdot \Lambda \cdot \left(
  \begin{array}{cc}
0 & \frac{1}{K_3} V_{\lambda^t}\\ \frac{1}{K_2} V_\lambda^{-1} & 0
  \end{array}\right)\,.
\]
Since $a \not= 0$ and $a \not= 1$, $V_\lambda = U_\lambda^t$,
$\lambda^t = \left[\frac{1}{a} : \frac{1}{b} : 1\right]$,
$U_{\lambda^t} = \frac{1}{a^2} \cdot U_\lambda$, $V_{\lambda^t} =
\frac{1}{a^2} U_\lambda^t$ (see Remark \ref{remark6.2}).

So $\Lambda^\prime = \left(
  \begin{smallmatrix}
    0 & -K_3a^2 U_\lambda^{-1}\\ -K_2U_\lambda & 0
  \end{smallmatrix}\right) \cdot \Lambda \cdot \left(
  \begin{smallmatrix}
    0 & \frac{1}{K_3} \cdot \frac{1}{a^2} \cdot U_\lambda^t\\ \frac{1}{K_2}
    (U_\lambda^{-1})^t & 0
  \end{smallmatrix}\right) = \left(
  \begin{smallmatrix}
    0 & k U_\lambda^{-1}\\ -\frac{1}{k}U_\lambda & 0
  \end{smallmatrix}\right) \cdot \Lambda \cdot \left(
  \begin{smallmatrix}
    0 & k U_\lambda^{-1} \\ -\frac{1}{k}U_\lambda & 0
  \end{smallmatrix}\right)^t$, where $k^2 = - a^2 \cdot \frac{K_3}{K_2}$.
\end{enumerate}
\end{proof}

In a similar way, we can prove the following lemma:




\begin{lem}\label{lemma6.8}
  Let $\lambda = [1:b:1] \in V(f_3)\smallsetminus \{P_0\}$.

  \begin{enumerate}
  \item Inside the family $\km_\lambda$, two matrices $\Lambda$ and
  $\Lambda^\prime$ are equivalent if and only if $\Lambda^\prime = T \cdot
  \Lambda \cdot T^t$, where
\[
T = \left(
  \begin{array}{cc}
K_4 \cdot \id & -K_3 \cdot U_\lambda^{-1}\\-K_2 \cdot U_\lambda &
K_1 \cdot \id
  \end{array}\right),\quad K_1, K_2, K_3, K_4 \in K \text{ such that } K_1 K_4
- K_2K_3 = 1\,.
\]

\item No $\lambda \in V(f_3)\smallsetminus \{P_0,[1:b:1]\}$ exists, such that
  a matrix from $\km_\lambda$ is equivalent to a matrix from $\km_{[1:b:1]}$.
  \end{enumerate}
\end{lem}




Now let us see ``how large'' the family $\km_\lambda$ is for a given $\lambda$
in $V(f_3)\smallsetminus\{P_0\}$.

For $\Lambda = \Lambda_{(\lambda,\Gamma)}$ in $\km_\lambda$, we
denote:
\[
\Gamma_1 = \left(
  \begin{array}{lll}
\phantom{-}0 & \phantom{-}a_1 & a_2\\
-a_1 & \phantom{-}0 & a_3\\
-a_2 & -a_3 & 0
  \end{array}\right),\quad \Gamma_2 = \left(
  \begin{array}{lll}
a_7 & a_8 & a_9\\
a_{10} & a_{11} & a_{12}\\
a_{13} & a_{14} & a_{15}
  \end{array}\right),\quad \Gamma_3 = \left(
  \begin{array}{lll}
\phantom{-}0 & \phantom{-}a_4 & a_5\\
-a_4 & \phantom{-}0 & a_6\\
-a_5 & -a_6 & 0
  \end{array}\right)\,.
\]
The condition $\det \Lambda = f^2$ provides 10 equations in the
above  15 parameters.   Six of these equations are linear in the
entries of $\Gamma_2$ and form a linear system of dimension three.

\begin{enumerate}
\item If $b = 0$ the solution of this system is:
\[
\begin{cases}
  a_7 = -a_{12} \cdot (a^2 + 1)\\
a_8 = a_{10} = a_{15} = 0\\
a_9 = a_{11} - a_{13}\\
a_{14} = a^2 \cdot a_{12}\,.
\end{cases}
\]

\item If $b \not= 0$, the system has the following solution:
\setlength{\extrarowheight}{5pt}
\[
\begin{cases}
  a_8 \; = -\frac{b}{a+1} a_7 + a_{15}\\
a_9 \; = - \frac{a-1}{b(a+1)} a_7 - \frac{a^2}{b^2} \cdot a_{15}\\
a_{10} = \frac{b}{a+1} \cdot a_7\\
a_{12} = -\frac{a^2+3}{(a+1)^2} \cdot a_7 + \frac{b}{a+1} a_{11} +
  \frac{1-a}{b(a+1)} a_{15}\\
a_{13} = \frac{a-1}{b(a+1)} \cdot a_7 + a_{11} + \frac{a^2}{b^2}
\cdot
a_{15}\\
a_{14} = \frac{2(1-a)}{(a+1)^2} \cdot a_7 - \frac{b}{a+1} \cdot
a_{11} +
  \frac{a-1}{b(a+1)} \cdot a_{15}\,.
\end{cases}
\]
\end{enumerate}

The other four equations are linear in the entries of $\Gamma_1$
with coefficients in $K[a_4, \dots, a_{15}]$ and have dimension
five:

\begin{alltt}
LIB"matrix.lib"; option(redSB);

ring r=0,(x(4),x(1),x(2),x(3),e,a,b,a(1..15)),dp; ideal
ii=a3+b3+1,e*b+a2-a+1,e*a+e-b2; qring s=std(ii);

matrix B[10][1]; B[1,1]=x(4)*a(1); B[2,1]=x(4)*a(2);
B[3,1]=-x(4)*a(7); B[4,1]=-x(4)*a(10)-(x(1)+x(3));
B[5,1]=x(4)*a(3); B[6,1]=-x(4)*a(8)-(x(1)-a*x(3));
B[7,1]=-x(4)*a(11)+x(2)+b*x(3); B[8,1]=-x(4)*a(9)-x(2)+b*x(3);
B[9,1]=-x(4)*a(12)+e*x(3); B[10,1]=x(4)*a(4);

matrix V[1][5]; V[1,1]=-x(4)*a(13)-x(2);
V[1,2]=-x(4)*a(14)-e*x(3); V[1,3]=-x(4)*a(15)+x(1)+(a-1)*x(3);
V[1,4]=x(4)*a(5); V[1,5]=x(4)*a(6);

poly p1=B[5,1]*B[10,1]-B[6,1]*B[9,1]+B[7,1]*B[8,1]; poly
p2=B[2,1]*B[10,1]-B[3,1]*B[9,1]+B[4,1]*B[8,1]; poly
p3=B[1,1]*B[10,1]-B[3,1]*B[7,1]+B[4,1]*B[6,1]; poly
p4=B[1,1]*B[9,1]-B[2,1]*B[7,1]+B[4,1]*B[5,1]; poly
p5=B[1,1]*B[8,1]-B[2,1]*B[6,1]+B[3,1]*B[5,1];

poly g=V[1,1]*p1-V[1,2]*p2+V[1,3]*p3-V[1,4]*p4+V[1,5]*p5; poly
f=x(4)^3+x(1)^3+x(2)^3+x(3)^3; g=g-f;

//For our skew symmetric matrix the condition g=f is equivalent to
//\(\det\Lambda = f\sp{2}\).

matrix H=coef(g,x(4)*x(1)*x(2)*x(3)); for(int j=1;j<=13;j++) {
H[1,j]=0; }

ideal I=H; I=interred(I);

I[1]=a(9)-a(11)+a(13) I[2]=a(8)+a(10)-a(15) I[3]=a(7)+a(12)+a(14)
I[4]=a*a(10)-e*a(11)+b*a(12)+2*e*a(13)+2*b*a(14)-2*a*a(15)+a(10)+a(15)
I[5]=2*e*a(10)+2*b*a(11)-2*a*a(12)-b*a(13)-a*a(14)-e*a(15)+a(12)+2*a(14)
I[6]=a(3)*a(4)-a(2)*a(5)+a(1)*a(6)+a(11)^2+a(10)*a(12)-a(11)*a(13)+a(13)^2
    -a(10)*a(14)-2*a(12)*a(15)-a(14)*a(15)
I[7]=a(1)*a(4)+a(3)*a(5)+a(2)*a(6)-a(10)^2+a(11)*a(12)+a(12)*a(13)+2*a(11)
    *a(14)-a(13)*a(14)+a(10)*a(15)-a(15)^2
I[8]=2*e^2*a(12)+2*a*b*a(12)-3*b^2*a(13)+2*e^2*a(14)-a*b*a(14)-3*e*b*a(15)
     -6*e*a(11)-b*a(12)+12*e*a(13)+2*b*a(14)-6*a*a(15)
I[9]=a(3)*a(5)*a(10)-a(2)*a(6)*a(10)-a(2)*a(5)*a(11)-a(1)*a(6)*a(11)+a(1)
     *a(5)*a(12)+a(3)*a(6)*a(12)-a(2)*a(5)*a(13)+2*a(1)*a(6)*a(13)+a(13)^3
     +a(2)*a(4)*a(14)+a(3)*a(6)*a(14)+a(10)*a(11)*a(14)+a(12)^2*a(14)-2
     *a(10)*a(13)*a(14)+a(12)*a(14)^2+a(3)*a(5)*a(15)+2*a(2)*a(6)*a(15)
     +a(11)*a(14)*a(15)-2*a(13)*a(14)*a(15)-a(15)^3-1
I[10]=2*e*a(2)*a(4)-2*e*a(1)*a(5)+2*b*a(2)*a(5)-2*a*a(3)*a(5)+2*b*a(1)*a(6)
     -4*a*a(2)*a(6)-2*b*a(11)^2+2*a*a(11)*a(12)+2*e*a(12)^2+5*b*a(11)*a(13)
     -4*a*a(12)*a(13)-2*b*a(13)^2-2*b*a(10)*a(14)-a*a(11)*a(14)+2*a*a(13)
     *a(14)-2*e*a(14)^2+3*e*a(11)*a(15)-6*e*a(13)*a(15)-2*b*a(14)*a(15)
     +6*a*a(15)^2+4*a(3)*a(5)+2*a(2)*a(6)-a(11)*a(12)+2*a(12)*a(13)+2*a(11)
     *a(14)-4*a(13)*a(14)-6*a(15)^2

ideal J=I[1],I[2],I[3],I[4],I[5],I[8];

//This is the ideal generated by the linear equations in the entries
//of \(\Gamma\sb{2}\).

ideal JJ=std(J); dim(JJ); 14

ideal J1=I[6],I[7],I[9],I[10];

//This is the ideal generated by the other four equations.

ideal JJ1=std(J1); dim(JJ1); 16
\end{alltt}

Let us summarize the results.

Let $M$ be an indecomposable graded rank 2, 6--generated MCM and
$\overline{M}$ the restriction of $M$ to the elliptic curve on our surface
defined by $f = x_4 = 0$.  Then $\overline{M} \cong N_\lambda \oplus
N_\lambda^\vee$ for a suitable 3--generated rank 1 MCM $N_\lambda = \coker
(\alpha_\lambda)$, $\lambda \in V(f,x_4) \smallsetminus \{[-1:0:1:0]\} \cong
V(f_3) \smallsetminus \{[-1:0:1]\} =: C$.  If $\lambda = [a:b:c]$ and
$\lambda^t := [c:b:a]$, then $N_\lambda^\vee \cong N_{\lambda^t}$, in
particular, there exist skew--symmetric $3 \times 3$--matrices $\Gamma_1,
\Gamma_3$ with constant entries not being zero simultaneously and a $3 \times
3$--matrix $\Gamma_2$ such that $M = \coker (\Lambda)$ for $\Lambda =
x_4\left(
  \begin{smallmatrix}
    \Gamma_1 & -\Gamma_2^t\\\Gamma_2 & \phantom{-}\Gamma_3
  \end{smallmatrix}\right) + \left(
  \begin{smallmatrix}
    0 & -\alpha_\lambda^t\\ \alpha_\lambda & \phantom{-}0
  \end{smallmatrix}\right)$, $\Gamma_1 = \left(
  \begin{smallmatrix}
    \phantom{-}0 & \phantom{-}a_1 & a_2\\-a_1 & \phantom{-}0 & a_3\\-a_2 &
    -a_3 & 0
  \end{smallmatrix}\right)$, $\Gamma_3 = \left(\begin{smallmatrix} \phantom{-}0 & \phantom{-}a_4
    & a_5\\-a_4 & \phantom{-}0 & a_6\\-a_5 & -a_6 &
    0\end{smallmatrix}\right)$, $\Gamma_2 = \left(
  \begin{smallmatrix}
    a_7 & a_8 & a_9\\a_{10} & a_{11} & a_{12}\\a_{13} & a_{14} & a_{15}
  \end{smallmatrix}\right)$ and $\det(\Lambda) = f^2$.

Let $\P = \P (2:2:2:-2:-2:-2)$ be the weighted projective space with respect
to the weights $2,2,2,-2,-2,-2$ and the coordinates
$(a_1:a_2:a_3:a_4:a_5:a_6)$.  Let $\A = \A^9$ be the 9--dimensional affine
space with the coordinates $(a_7, \dots, a_{15})$.

A point $(\lambda; \underline{a}) \in C \times \P \times \A$ corresponds to an
equivalence class of matrices $\Lambda = x_4 \left(
  \begin{smallmatrix}
    \Gamma_1 & -\Gamma_2^t\\\Gamma_2 & \phantom{-}\Gamma_3
  \end{smallmatrix}\right) + \left(
  \begin{smallmatrix}
    0 & -\alpha_\lambda^t\\\alpha_\lambda & \phantom{-}0
  \end{smallmatrix}\right)$ under the action of the group $\{U_k \mid k \in
K^\ast\}$, $\Lambda \mapsto U_k \Lambda U_k^t$.   The duality of $C$ defined
by the 3--generated MCM of rank 1, $\lambda = [a:b:c] \mapsto [c:b:a] =
\lambda^t$, induces an $S_2$--action on $C \times \P \times \A$.   In terms of
matrices it is defined by $\Lambda \mapsto U \Lambda U^t$, $U = \left(
  \begin{smallmatrix}
    \phantom{-}0 & U_\lambda^{-1}\\-U_\lambda & 0
  \end{smallmatrix}\right)$.   Let $\km \subseteq C \times \P \times \A$ be
the $S_2$--invariant closed subset defined by $\det(\Lambda) = f^2$.   Let
$\pi : \km \to C$ be the canonical projection.

\begin{thm}\label{theorem6.10}

  \begin{enumerate}
  \item Every indecomposable graded rank 2, 6--generated MCM is represented
  by a point in $\km$.

\item $\km \smallsetminus \pi^{-1} (\{[1:b:1] \mid b^3 = -2\})/S_2$ is the
  moduli space  of isomorphism classes of indecomposable graded rank 2,
  6--generated MCM $M$ such that the restriction to $V(f,x_4)$, $\overline{M}
  \cong N_\lambda \oplus N_\lambda^\vee$ for $N_\lambda$ being not self--dual.
  This moduli space is 5--dimensional.

\item Let $H = \left\{\left(
      \begin{smallmatrix}
        \phantom{-}K_4 \id & & -K_3 U_\lambda^{-1}\\
-K_2U_\lambda & & \phantom{-}K_1 \id
      \end{smallmatrix}\right),\; K_1, K_2, K_3, K_4 \in K,\; K_1 K_4
    -K_2K_3=1 \right\}$, then $H$ acts on $\pi^{-1}(\{[1:b:1] \mid b^3 =-2\})$
     and $\pi^{-1} (\{[1:b:1] \mid b^3 = -2\})/H$ is the moduli space of
     isomorphism classes of indecomposable graded rank 2, 6--generated MCM $M$
     such that the restriction to $V(f,x_4)$, $\overline{M} \cong N_\lambda \oplus
     N_\lambda$ for $N_\lambda$ being self--dual.   This moduli space is
     2--dimensional.
  \end{enumerate}
\end{thm}

\begin{rem}\label{rem6.12} It is well known that the ideal  defining
  5 general points in ${\bf P}_K^3$ is Gorenstein (this means any four from
  them are not on a hyperplane). Restricting to the 5 general points on the
  surface $V(f)$ we get a family of Gorenstein ideals whose isomorphism
  classes of 2-syzygies over $R$ (they are indecomposable, graded, rank two,
  6-generated MCM modules) form a 5-parameter family (see \cite{Mi},
  \cite{IK}). Here we give an example. Let \mbox{$[1:0:0:-1]$},
  \mbox{$[1:0:-1:0]$}, \mbox{$[1:-1:0:0]$}, \mbox{$[1:-u:0:0]$},
  \mbox{$[1:-u:1:-u]$}, \mbox{$u^2+u+1=0$} be 5 general points on $V(f)$ and
  $I$ the ideal defined by these points in $R$. $I$ is generated by the
  following quadratic forms: $x_2x_4+ ux_3x_4$, $-ux_2x_3+ux_3x_4$,
  $x_1x_4+x_4^2-(1-u)x_3x_4$, $u(x_1+x_3)x_3 +2x_3x_4$,
  $-x_3x_4-x_1^2+ux_1x_2-u^2x_2^2+x_3^2+x_4^2$.  Then the second syzygy of $I$
  over $R$ is the cokernel of a skew symmetric matrix $A$ defined by

$A[1,1]=A[2,2]=A[3,3]=A[4,4]=A[5,5]=A[6,6]=0$,\\
$A[1,2]=(-3u-2)x_3+(2u-1)x_4=-A[2,1]$,\\
$A[1,3]=-ux_1+(-2u+1)x_2+(u+1)x_3+ux_4=-A[3,1]$,\\
$A[1,4]=(u-2)x_1-x_2+(-3u-4)x_3+(2u-1)x_4=-A[4,1]$,\\
$A[1,5]=(u+1)x_3-ux_4=-A[5,1]$,\\
$A[1,6]=-ux_1+(u+1)x_2+(1/7u+3/7)x_3+(-3/7u-2/7)x_4=-A[6,1]$,\\
$A[2,3]=(u-2)x_1-x_2+x_3+(-u+2)x_4=-A[3,2]$,\\
$A[2,4]=(3u+2)x_1+(2u+3)x_2+4ux_3+x_4=-A[4,2]$,\\
$A[2,5]=(-3u-1)x_3+(u-2)x_4=-A[5,2]$,\\
$A[2,6]=(-u-2)x_1+(-u+1)x_2+(-u-1)x_3+ux_4=-A[6,2]$,\\
$A[3,4]=-3x_3=-A[4,3]$,\\
$A[3,5]=(u+1)x_3=-A[5,3]$, \\
$A[3,6]=(-6/7u-4/7)x_3+x_4=-A[6,3]$,\\
$A[4,5]=(-3u-1)x_3=-A[5,4]$,\\
$A[4,6]=-ux_3+ux_4=-A[6,4]$,\\
$A[5,6]=-x_1-ux_2=-A[6,5]$.
\end{rem}

\end{document}